\documentclass[english,a4paper]{amsart}

\usepackage[english,lined,linesnumbered,ruled]{algorithm2e}
\usepackage{amsfonts,bm,bbm} 
\usepackage{amsmath,amsopn,amssymb,amsthm,mathrsfs,mathtools}
\usepackage{array}
\usepackage{caption}
\usepackage{cleveref} 
\usepackage{dsfont} 
\usepackage{enumitem} 
\usepackage{multirow}
\usepackage{stmaryrd} 
\usepackage{tikz}
\def \ind{\mathds{1}}

\def \CV{\mathrm{CV}}
\def \EE{\mathbb{E}}
\def \LL{\mathbb{L}}
\def \NN{\mathbb{N}}
\def \RR{\mathbb{R}}
\def \PP{\mathbb{P}}

\def \var{\mathrm{Var}}

\def \cA{\mathcal{A}}
\def \cB{\mathcal{B}}

\def \cD{\mathcal{D}}

\def \cI{\mathcal{I}}

\def \cM{\mathcal{M}}
\def \cN{\mathcal{N}}
\def \cO{\mathcal{O}}

\def \cR{\mathcal{R}}
\def \cS{\mathcal{S}}
\def \cU{\mathcal U}
\def \cX{\mathcal X}

\def \bx{\mathbf x}
\def \bX{\mathbf X}

\def \bz{\mathbf z}

\newcommand{\eps}{\varepsilon}

\newtheorem{hyp}{Assumption}[section] 
\newtheorem{theo}{Theorem}[section]
\newtheorem*{lem*}{Lemma}
\newtheorem{lem}[theo]{Lemma}
\newtheorem{mydef}{Definition}[section]
\newtheorem{prop}[theo]{Proposition}
\newtheorem{remark}{Remark}[section]

\crefname{hyp}{assumption}{assumptions}
\Crefname{hyp}{Assumption}{Assumptions}
\crefname{lem}{lemma}{lemmas}
\Crefname{lem}{Lemma}{Lemmas}
\crefname{prop}{proposition}{propositions}
\Crefname{prop}{Proposition}{Propositions}
\crefname{theo}{theorem}{theorems}
\Crefname{theo}{Theorem}{Theorems}

\numberwithin{equation}{section}

\newenvironment{myproof}[1][\proofname]{%
  \begin{proof}[#1]$ $\par\nobreak\ignorespaces
}{%
  \end{proof}
}

\newcommand{\vero}[1]{{\color{black}#1}}
\newcommand{\kev}[1]{{\color{black}#1}}
\newcommand{\verob}[1]{{\color{black}#1}}
\newcommand{\veroc}[1]{{\color{black}#1}}

\begin{document}


\title[Consistency of Conditional Distribution Forest]{Random Forest Estimation of Conditional Distribution Functions and Conditional Quantiles}


\author{Kévin Elie-Dit-Cosaque }
\address{SCOR France. }
\author{Véronique Maume-Deschamps}
\address{Univ Lyon, Université Claude Bernard Lyon 1,\\ CNRS UMR 5208, Institut Camille Jordan, \\
F-69622 Villeurbanne, France}



\begin{abstract}
We propose a theoretical study of two realistic estimators of \veroc{conditional distribution functions using random forests}. The estimation process uses the bootstrap samples generated from the original dataset when constructing the forest. Bootstrap samples are reused to define the first estimator, while the second requires only the original sample, once the forest has been built. We prove that both proposed estimators of the conditional distribution functions are consistent uniformly a.s. \vero{To the best of our knowledge, it is the first proof of a.s. consistency and including the bootstrap part. The consistency result holds for a large class of functions, including additive models and products.} \veroc{The consistency of conditional quantiles estimators follows that of distribution functions estimators using standard arguments. } 
\end{abstract}
%






\maketitle

\section{Introduction}

Conditional distribution functions and conditional quantiles estimation is an important task in several domains including environment, insurance or industry. It is also an important tool for Quantile Oriented Sensitivity Analysis (QOSA), see e.g., \cite{fort2016new,maume2018estimation,browne2017estimate,qosa_shapley}. In order to estimate conditional quantiles, various methods exist such as kernel based estimation or quantile regression \cite{koenker2001quantile} but they present some limitations. Indeed, the performance of kernel methods strongly depends on the bandwidth parameter selection and quickly breaks down as the number of covariates increases. On the other hand, quantile regression \vero{as presented in \cite{koenker2001quantile}  may not be adapted} in a non-gaussian setting since the true conditional quantile  is not necessarily a linear combination of the input variables \cite{maume2017quantile}. To overcome these issues, we propose to explore the random forest estimation of conditional quantiles \cite{meinshausen2006quantile}.

Random forest algorithms allow a flexible modeling of interactions in high dimension by building a large number of regression trees and averaging their predictions. The most famous random forest algorithm is that of \cite{breiman2001random} whose construction is based on the seminal work of \cite{amit1997shape,ho1998random,dietterich2000ensemble}. Breiman's random forest estimate is a combination of two essential components: Bagging and Classification And Regression Trees (CART)-split criterion \cite{breiman1984classification}. Bagging for \textit{bootstrap-aggregating} was proposed by \cite{breiman1996bagging} in order to improve the performance of weak or unstable learners. \vero{Two types of randomness are included: a \kev{bootstrap} sample is used to construct each tree and at each split step, \verob{some input variables are chosen randomly in $\{X_1\/,\ldots\/, X_d\}$. }}

Random forests are also related to some local averaging algorithms such as nearest neighbors methods \cite{lin2006random,biau2010layered} or kernel estimators \cite{scornet2016random}. More precisely, thanks to \cite{lin2006random}, the random forest method can be seen as an adaptive neighborhood regression procedure and therefore the prediction (estimation of the conditional mean) can be formulated as a weighted average of the observed response variables. 

Based on that approach, we develop a Weighted Conditional Empirical Cumulative Distribution Function (W\textunderscore C\textunderscore ECDF) approximating the Conditional Cumulative Distribution Function (C\textunderscore CDF). Then, $\alpha$-quantile estimates are obtained by using W\textunderscore C\textunderscore ECDF instead of C\textunderscore CDF. \cite{meinshausen2006quantile} defined a W\textunderscore C\textunderscore ECDF with weights using the original dataset whereas we allow to construct the weights using the bootstrap samples, as it is done practically in regression random forests. \vero{We prove the almost sure consistency of these \veroc{two} estimators. The main hypothesis are:
\begin{enumerate}
 \item the convergence to $0$ of the variation of the conditional distribution function on leaves;
 \item a control of the size of the leaves, implying that the trees are not fully developed.
\end{enumerate}
In Section \ref{sec:variation}, we prove that the first hypothesis is satisfied for a large class of models, including additive models and product functions. \veroc{Also, this class of functions is dense in the set of non-negative continuous functions on $[0\/,1]^d$.} Moreover, the first hypothesis \veroc{is  satisfied} for some modified version of the CART algorithm, such as in \cite{meinshausen2006quantile}. The second hypothesis may be seen as a stopping rule in the algorithm, so that it should be controlled by practitioners. }
\vero{An} implementation of both algorithms is made available within a \texttt{Julia} package called \texttt{ConditionalDistributionForest} \cite{ConditionalDistributionForest} as well as a \texttt{Python} package named \texttt{qosa-indices}, \cite{qosa_indices}.

The C\textunderscore CDF can be seen as a regression function. 
Several authors such as \cite{breiman2004consistency, biau2012analysis, wager2015adaptive, scornet2015consistency, mentch2016quantifying, wager2018estimation,goehry2019random} have established asymptotic properties of particular variants and simplifications of the original Breiman's random forest algorithm. \veroc{For instance, in \cite{biau2012analysis}, the tree construction is done by choosing at random the splitting coordinate and once the latter is selected, the split is at the midpoint of the chosen side.}  Facing some theoretical issues with the bootstrap, most studies replace it by subsampling, assuming that each tree is grown with $s_n<n$ observations randomly chosen without replacement from the original dataset. Most of the time, in order to ensure the convergence of the simplified model, the subsampling rate $s_n / n$ is assumed to tend to zero at some prescribed rate, assumption that excludes the bootstrap mode. Besides, consistency is generally showed by assuming that the number of trees goes to infinity which is not fully relevant in practice. Under some conditions, \cite{scornet2016asymptotics} showed that if the infinite random forest regression estimator is $\LL^2$ consistent then so does the finite random forest regression estimator when the number of trees goes to infinity in a controlled way. \\
Recent attempts to bridge the gap between theory and practice, provide some results on random forest algorithms at the price of fairly strong conditions. For example, \cite{scornet2015consistency} showed the $\LL^2$ consistency of random forests in an additive regression framework by replacing the bootstrap step by subsampling. Their result rests on a fundamental lemma developed in \cite{scornet2015supplementary} which reviews theoretical random forest, \vero{the additive assumption being required here}. 
Furthermore, consistency and asymptotic normality of the whole algorithm were recently proved under strong conditions by \cite{wager2018estimation} replacing bootstrap by subsampling and simplifying the splitting step. One of the strong conditions used in Theorem 3.1. of \cite{wager2018estimation} is that the individual trees satisfy a condition called \textit{honesty}. An example of an honest tree given by the authors is \textit{one where the tree is grown using one subsample, while the predictions at the leaves of the tree are estimated using a different subsample.} Due to this assumption, the authors admit that their theorems are not valid for the practical applications most of the time because almost all implementations of random forests use the training sample twice. \veroc{In a recent work \cite{klusowki}, the $\LL^2$ consistency in a regression framework is also obtained for the expectation over all possible trees of the regression function (infinite forest),  the bootstrap process is not taken into account. The result in \cite{klusowki} holds for a class of regression functions broader than the additive model. This class of functions requires some control on the partial derivatives of the model.}

Thus, despite an active investigation during the last decade, \vero{further consistency results are still welcome. }
Our major contribution is the proof of the almost sure uniform convergence of the estimator W\textunderscore C\textunderscore ECDF both using the bootstrap samples (\Cref{theo:4:num_1}) or the original one (\Cref{theo:4:num_2}).  \vero{Remark} that \cite{meinshausen2006quantile} gave a proof of the consistency in probability of the W\textunderscore C\textunderscore ECDF for a simplified model where the weights are considered as constant while they are indeed random variables heavily data-dependent. \vero{We provide an a.s. consistency proof under realistic assumptions for a method based on bootstrap samples. This consistency result holds for a large class of functions called $\spadesuit$-{\bf class} which contains, for instance, additive functions, product functions. \veroc{Our $\spadesuit$-class also contains some dense class of functions, as does the class considered in \cite{klusowki}. Contrary to \cite{klusowki}, our result holds for finite random forests and takes the bootstrap into account, so that it is closer to practical situations.}  Also, a sub-product of the consistency proof is an asymptotic proximity result between empirical and theoretical trees (see Proposition \ref{prop:close}), which is interesting in itself for further theoretical studies on random forests. }

The paper is organized as follows. Breiman's random forest algorithm is detailed in Section \ref{sec:2:RF} and notations are stated. The random forest estimations of C\textunderscore CDF based both on bootstrap samples and the original dataset are introduced in Section \ref{sec:3:QRF} as a natural generalization of regression random forests. The main consistency results are presented in Section \ref{sec:4:consistency} \vero{and a deep discussion on the variation of the conditional distribution function on leaves is proposed in Section \ref{sec:variation}. The main proofs are gathered in Section \ref{sec:5:proofs} \kev{and some of them are postponed in Appendix \ref{sec:proof2}}. A short conclusion is given in Section \ref{sec:7:conclusion}.}

\section{Breiman's random forest}\label{sec:2:RF}

The aim of this section is to present  Breiman's random forest algorithm as well as notations used throughout this paper. 

Random forest is a generic term to name an aggregation scheme of decision trees allowing to deal with both supervised classification and regression tasks and we focus on the latter in this paper. \\
The general framework is the nonparametric regression estimation where an input random vector \kev{$\bX \in \cX = \displaystyle\prod_{i=1}^d [u_i\/,v_i] \subset \RR^d$} is observed and a response $Y \in \RR $ is predicted by estimating the regression function $m(\bx) = \EE \left[ \left. Y \right| \bX = \bx \right]$. 
We assume that we are given a training sample $\cD_n = \left( \bX^j,Y^j \right)_{j=1,\ldots,n}$ of independent random variables distributed as the prototype pair $\left( \bX,Y\right)$ which is a $(d+1)$-dimensional random vector. The purpose is to use the dataset $\cD_n$ to construct an estimator $m_n:\cX \to \RR$ of the function $m$.

Random forests proposed by \cite{breiman2001random} build a predictor consisting of a collection of $k$ randomized regression trees grown based on the CART algorithm. \\
The CART-split criterion of \cite{breiman1984classification} is used in the construction of the individual trees to recursively partition the input space $\cX$ in a dyadic manner. More precisely, at each step of the partitioning, a part of the space is divided into two sub-parts according to the best cut perpendicular to the axes. This best cut is selected in each node of the tree by optimizing the CART-split criterion over the $d$ variables, i.e. minimizing the prediction squared error in the two child nodes.  The trees are \vero{then grown recursively} until reaching a stopping rule. There are several rules, but one generally proposed is that the tree construction continues while leaves contain at least $min\mathunderscore samples\mathunderscore leaf$ elements. This criterion is implemented in the \texttt{RandomForestRegressor} class of the \texttt{python} package \texttt{Scikit-Learn} \cite{scikit-learn} or in the \texttt{build\textunderscore forest} function of the \texttt{Julia} \cite{bezanson2017julia} package \texttt{DecisionTree}.

Building several different trees from a single dataset requires to randomize the tree building process. \cite{breiman2001random} proposed to inject some randomness both in the dataset and in the tree construction. First of all, prior to the construction of each tree, a resampling step is done by bootstrapping \cite{Efron_1979} from the original dataset, that is, by randomly sampling $n$ observations with replacement. Only these bootstrap observations are taken into account in the tree building. Accordingly, the $min\mathunderscore samples\mathunderscore leaf$ hyperparameter introduced previously refers, in the random forest method, to the minimum number of bootstrap observations contained in each leaf of a tree. Secondly, at each step of the tree construction, instead of optimizing the CART-split criterion over the $d$ variables, \veroc{a number of variables called $max\mathunderscore features$ is selected uniformly at random among the $d$ variables, the set of selected variables is denoted $\cM_{\mbox{try}}$. In what follows, we shall assume that the $max\mathunderscore features$ co-variates are selected randomly, with positive probability for each co-variate to be selected, which includes of course the uniform selection. } Then, the best split is chosen as the one optimizing the CART-split criterion \vero{as follows: for a given node $A= \displaystyle\prod_{i=1}^d [a_i\/,b_i]$,  the CART-split is given by maximizing} \kev{over $j\in \cM_{\mbox{try}}$ and $z\in A^j= [a_j\/,b_j]$
\begin{eqnarray}\label{def:L_n}
L_A^n(j,z) & = &\dfrac1{N_n^b(A)} \sum_{i=1}^n \left( Y^i - \overline{Y}_A \right)^2 \ind_{\left\lbrace \bX^i \in A \right\rbrace}-   \\
&&\dfrac1{N_n^b(A)} \sum_{i=1}^n \left( Y^i-\overline{Y}_{A_L}\ind_{\left\lbrace X_j^i \leq z \right\rbrace} - \overline{Y}_{A_R} \ind_{\left\lbrace X_j^i > z \right\rbrace} \right)^2  \ind_{\left\lbrace \bX^i \in A \right\rbrace}\/,\nonumber
\end{eqnarray}
where $A_L = \left\lbrace \bx \in A, x_j \leq z \right\rbrace$ (resp. $A_R$) is the left (resp. right) child of $A$, $\overline{Y}_B$ is the mean of the $Y^i$'s in the bootstrap sample with $\bX^i \in B$} and $N_n^b(B)$ denotes the number of elements in the bootstrap sample belonging to $B$.\\ 
\vero{Now, for} any query point $\bx \in \cX$, the $\ell$-th tree estimates $m(\bx)$ as follows
\begin{equation}\label{eq:2:EC_tree_estimator}
m_n^b \left( \bx;\Theta_\ell,\cD_n \right) = \sum_{j \in \cD_n^{\star} \left( \Theta_\ell \right)} \dfrac{\ind_{\left\lbrace \bX^j \in A_n \left( \bx;\Theta_\ell,\cD_n \right) \right\rbrace}}{N_n^b \left( \bx;\Theta_\ell,\cD_n \right)} Y^j \ ,
\end{equation}
where:
\begin{itemize}
\item $ \Theta_{\ell} ,\ell=1,\ldots,k$ are independent random vectors, distributed as a generic random vector $\Theta = \left( \Theta^1, \Theta^2 \right)$ and independent of $\cD_n$. $\Theta^1$ contains indexes of observations that are used to build each tree, i.e. the bootstrap sample and $\Theta^2$ indexes of splitting candidate variables in each node, \veroc{we assume that $\Theta^2$ gives a positive probability to each co-variate,}
\item $\cD_n^{\star} \left( \Theta_\ell \right)$ is the bootstrap sample selected prior to the tree construction,
\item $A_n \left( \bx;\Theta_\ell,\cD_n \right)$ is the tree cell (subspace of $\cX$) containing $\bx$,
\item $N_n^b \left( \bx;\Theta_\ell,\cD_n \right)$ is the number of elements of $\cD_n^{\star} \left( \Theta_\ell \right)$ that fall into $A_n \left( \bx;\Theta_\ell,\cD_n \right)$. 
\end{itemize}
The trees are then combined to form the finite forest estimator
\begin{equation}\label{eq:2:EC_forest_estimator}
m_{k,n}^b \left( \bx;\Theta_1,\ldots,\Theta_k,\cD_n \right) = \dfrac{1}{k}\sum_{\ell=1}^{k} m_n^b \left( \bx;\Theta_\ell,\cD_n \right) \ .
\end{equation}

We may now present the conditional distribution function estimators.

\section{Conditional Distribution Forests}\label{sec:3:QRF}

We aim to estimate $F(y|\bx) = \PP(Y \leqslant y | \bX=\bx)$. Two estimators may be defined. One uses the bootstrap samples both in the forest construction and in the estimation. The other uses 
the original sample in the estimation part. Once the distribution function has been estimated, the conditional quantiles may be estimated straightforwardly. 

\subsection{Bootstrap samples based estimator}\label{sub_sec:3.1:QRF}
\veroc{First of all, let us define the random variable $B_{j} \left( \Theta_\ell^1,\cD_n \right)$ as the number of times that the observation $\left( \bX^j, Y^j \right)$ has been drawn from the original dataset in the bootstrap sample for the $\ell$-th tree construction. 
The conditional mean estimator in Equation \eqref{eq:2:EC_forest_estimator} rewrites as}
\begin{eqnarray*}\label{eq:3:bootsrap_weighted_EC_forest_estimator}
m_{k,n}^b \left( \bx;\Theta_1,\ldots,\Theta_k,\cD_n \right) & = &\sum_{j=1}^{n} \left( \dfrac{1}{k} \sum_{\ell=1}^{k}  \frac{B_{j} \left( \Theta_\ell^1,\cD_n \right)\ind_{\left\lbrace \bX^j \in A_n \left( \bx;\Theta_\ell,\cD_n \right) \right\rbrace}}{ N_n^b \left( \bx;\Theta_\ell,\cD_n \right)} \right) Y^j \nonumber \\
& = &\sum_{j=1}^{n} w_{n,j}^b \left( \bx; \Theta_1,\ldots,\Theta_k,\cD_n \right) Y^j \ ,
\end{eqnarray*}
where the weights are defined by
\veroc{\begin{equation}\label{eq:ch4:3:bootstrap_weights}
w_{n,j}^b \left( \bx; \Theta_1,\ldots,\Theta_k,\cD_n \right) = \dfrac{1}{k} \sum_{\ell=1}^{k} \frac{B_{j} \left( \Theta_\ell^1,\cD_n \right)\ind_{\left\lbrace \bX^{j} \in A_n \left( \bx;\Theta_\ell,\cD_n \right) \right\rbrace}}{N_n^b \left( \bx;\Theta_\ell,\cD_n \right)} \ .
\end{equation}
}%
Note that the weights $w_{n,j}^b \left( \bx; \Theta_1,\ldots,\Theta_k,\cD_n \right)$ are non-negative random variables as functions of $\Theta_1,\ldots,\Theta_k, \cD_n$ and their sum for 
$j=1, \ldots, n$ equals $1$.

The random forest estimator \eqref{eq:3:bootsrap_weighted_EC_forest_estimator} can be seen as a local averaging estimate. Indeed, as mentioned by \cite{scornetpromenade}, the regression trees make an average of the observations located in a neighborhood of $\bx$, this neighborhood being defined as the leaf of the tree containing $\bx$. The forest, which aggregates several trees, also operates by calculating a weighted average of the observations in a neighborhood of $\bx$. However, in the case of forests, this neighborhood results from the superposition of the neighborhoods of each tree, and therefore has a more complex shape. Several works have tried to study the random forest algorithm from this point of view (local averaging estimate) such as \cite{lin2006random} who was the first to point out the connection between the random forest and the adaptive nearest-neighbors methods, further developed by \cite{biau2010layered}. Some works such as \cite{scornet2016random} have also studied random forests through their link with the kernel methods.

We are interested in the Conditional Cumulative Distribution Function (C\textunderscore CDF) of $Y$ given $\bX = \bx$ in order to obtain the conditional quantiles. Pairing the following equality
\begin{equation}\label{eq:3:CDF_as_regression_function}
F \left( \left. y \right| \bX=\bx \right) = \PP \left( \left. Y \leqslant y \right| \bX = \bx \right) = \EE \left[ \left. \mathds{1}_{ \{Y \leqslant y \}} \right| \bX = \bx \right] \ ,
\end{equation}
with the weighted approach described above, we propose to estimate the C\textunderscore CDF as follows
\begin{equation}\label{eq:3:bootsrap_CDF_estimator}
F_{k,n}^b \left( \left. y \right| \bX=\bx;\Theta_1,\ldots,\Theta_k,\cD_n \right) = \sum_{j=1}^{n} w_{n,j}^b \left( \bx; \Theta_1,\ldots,\Theta_k,\cD_n \right) \ind_{ \{Y^j \leqslant y\}} \ .
\end{equation}
Hence, given a level $\alpha \in \left] 0,1 \right[$, the conditional quantile estimator $\widehat{q}^\alpha \left( \left. Y \right| \bX = \bx \right)$ is defined as follows
\[
\widehat{q}^\alpha \left( \left. Y \right| \bX = \bx \right) = \inf \left\lbrace Y^{p}, p=1,\ldots,n : F_{k,n}^b \left( \left. Y^p \right| \bX=\bx;\Theta_1,\ldots,\Theta_k,\cD_n \right) \geqslant \alpha \right\rbrace \ .
\]

Let us turn now to the estimator using the original sample.

\subsection{Original sample based estimator}\label{sub_sec:3.2:QRF}

Trees are still grown with their respective bootstrap sample $\cD_n^{\star} \left( \Theta_\ell \right), \ell = 1, \ldots, k$. But instead of considering them in the estimation, we may use the original sample $\cD_n$. Consider the weights
\begin{equation}\label{eq:ch4:3:original_weights}
w_{n,j}^o \left( \bx; \Theta_1,\ldots,\Theta_k,\cD_n \right) = \dfrac{1}{k} \sum_{\ell=1}^{k} \dfrac{\ind_{\left\lbrace \bX^j \in A_n \left( \bx;\Theta_\ell,\cD_n \right) \right\rbrace}}{N_n^o \left( \bx;\Theta_\ell,\cD_n \right)} \ ,
\end{equation}
where $N_n^o \left( \bx;\Theta_\ell,\cD_n \right)$ is the number of points of $\cD_n$ that fall into $A_n \left( \bx;\Theta_\ell,\cD_n \right)$. As previously, the weights $w_{n,j}^o \left( \bx; \Theta_1,\ldots,\Theta_k,\cD_n \right)$ are non-negative random variables as functions of $\Theta_1,\ldots,\Theta_k, \cD_n$ and their sum over $j=1, \ldots, n$ equals 1.

It was proposed in \cite{meinshausen2006quantile} to estimate the C\textunderscore CDF with
\begin{equation}\label{eq:3:CDF_estimator}
F_{k,n}^o \left( \left. y \right| \bX=\bx;\Theta_1,\ldots,\Theta_k,\cD_n \right) = \sum_{j=1}^{n} w_{n,j}^o \left( \bx; \Theta_1,\ldots,\Theta_k,\cD_n \right) \ind_{ \{Y^j \leqslant y\}} \ .
\end{equation}
The conditional quantiles are then estimated by plugging $F_{k,n}^o \left( \left. y \right| \bX=\bx;\Theta_1,\ldots,\Theta_k,\cD_n \right)$ instead of $F \left( \left. y \right| \bX=\bx \right)$ as before.
A \texttt{python} library named \texttt{qosa-indices} has also been developed to perform the numerical estimations of conditional distributions and quantiles for both methods. It is available at \cite{qosa_indices} and uses \texttt{Scikit-Learn, Numpy, Numba}. Both approaches are also implemented in a \texttt{Julia} package based on the library \texttt{DecisionTree} and that is available at \cite{ConditionalDistributionForest}.

It has to be noted that a package called \texttt{quantregForest} has been made available in \texttt{R} \cite{Rsoftware} and can be found at \cite{qrf}. The estimation method currently implemented in \texttt{quantregForest} is different from the method described in \cite{meinshausen2006quantile}. As a matter of fact, for a new observation $\bx$ and the $\ell$-th tree, one element of $\cD_n = \left( \bX^j, Y^j \right)_{j=1,\ldots,n}$ falling into in the leaf node $A_n \left( \bx;\Theta_\ell,\cD_n \right)$ is chosen at random. This gives, $k$ values of $Y$ and allows to estimate the conditional distribution function with the classical Empirical Cumulative Distribution Function associated with the empirical measure. 

\section{Consistency results}\label{sec:4:consistency}

In this section, we state our main results, which are the uniform a.s. consistency of both estimators $F_{k,n}^b$ and $F_{k,n}^o$ of the conditional distribution function. It constitutes the most interesting result of this paper because it handles the bootstrap component and gives the almost sure uniform convergence. Indeed, most of the studies \cite{scornet2015consistency,wager2018estimation,goehry2019random} replaces the bootstrap by subsampling without replacement in order to avoid the mathematical difficulties induced by this one and therefore differs slightly from the procedure used in practice. \veroc{In \cite{klusowki}, the bootstrap procedure is not analyzed and the $\LL^2$ consistency is proved for the expectation over $\Theta$ of the estimator of the regression function, which does not lead directly to the result for finite forests.}

\cite{meinshausen2006quantile} showed the uniform convergence in probability of a simplified version of the estimator $F_{k,n}^o$. In \cite{meinshausen2006quantile}, the weights $w_{n,j}^o \left( \bx; \Theta_1,\ldots,\Theta_k,\cD_n \right)$ are in fact considered to be non-random while they are indeed random variables depending on $\left( \Theta_\ell \right)_{\ell=1,\ldots,k}$ and $\cD_n$.

Overall, proving the consistency of the forest methods whose construction depends both on the $\bX^j$’s and on the $Y^j$’s is a difficult task. This feature makes the resulting estimate highly data-dependent, and therefore difficult to analyze. A simplification widely used by most authors from a theoretical point of view is to work with random forest estimates whose tree shape depends only on $\bX^j$’s which \cite{devroye2013probabilistic} called the $\bX$-property but the $Y^j$’s are still used to compute the prediction, either the conditional mean or the conditional distribution function for example. One of the first results dealing with data-dependent random forest estimator of the regression function is \cite{scornet2015consistency} who showed the $\LL^2$ consistency in an additive regression framework by replacing the bootstrap by subsampling. \veroc{Our main result is stated in \Cref{theo:4:num_1,theo:4:num_2} hereafter but we first rewrite it in a more explicit way in order to emphasize their applicability. We begin by a description of the class of functions on which our result holds, then we state the hypothesis and the result. \\
In what follows, $\cX$ is a compact hyper-rectangle of $\RR^d$: \\
$\cX = \displaystyle \prod_{i=1}^d [u_i\/,v_i]$, $-\infty < u_i \leq v_i <\infty$ and we denote by $\cA$ the set of hyper-rectangles in $\cX$: $A\in \cA$ writes $A=\displaystyle\prod_{i=1}^d [a_i\/,b_i]$ with $u_i\leq a_i \leq b_i \leq v_i$. Also,  we denote by $A^{-j} = \displaystyle \prod_{k\neq j} [a_k, b_k]$ and $A^{J} = \displaystyle\prod_{k\in J} [a_k\/,b_k]$ for any $J \subset \{1\/, \ldots \/, d\}$.  Given $\bx \in \RR^d$, $\bx^{-j}$ is the vector of $\RR^{d-1}$ where the $j$-th coordinate has been removed and $\bx^J$ is the vector of $\RR^J$ whose coordinates are $x_j$, $j\in J$. 
\begin{mydef}\label{def:pique}
 Let $f:\cX \to \RR$, it does NOT belong to the  $\spadesuit$-class if there exists a rectangle $A=\displaystyle\prod_{j=1}^d [a_j\/,b_j] \subset \cX$, with $a_j\leq b_j$ such that for all $j=1\/,\ldots\/,d$,  $z \mapsto \EE \left[ f \left( z, \bX^{-j} \right)\ind_{\left\lbrace \bX^{-j} \in A^{-j} \right\rbrace} \right]$ is constant on $[a_j,b_j]$ and \kev{$f$ is not constant on $A$}.
\end{mydef}
\begin{remark}
From Definition \ref{def:pique}, if $f$ belongs to the $\spadesuit$-class then either  for  any rectangle $A$, there exists $j=1\/,\ldots\/, d$ such that $z\mapsto \EE \left[ f \left( z,\bX^{-j} \right) \ind_{\left\lbrace \bX^{-j}\in A^{-j} \right\rbrace} \right]$ is not constant on $[a_j,b_j]$ or if for some rectangle $A$, $z\mapsto \EE \left[ f \left(z,\bX^{-j} \right)\ind_{\left\lbrace \bX^{-j} \in A^{-j} \right\rbrace} \right]$ is constant for all $j=1\/,\ldots\/, d$ then $f$ is constant on $A$. 
\end{remark}
Let us give some examples of functions of the $\spadesuit$-class, the proof is more or less straightforward, a sketch of proof is given in the Appendix (see Lemma \ref{lem:pique}). \\
\begin{remark}\label{rem:pique}The $\spadesuit$-class contains:
\begin{enumerate}
 \item additive functions $f(\bx) = \sum_{j=1}^d f_j(x_j)$,
 \item product functions $f(\bx) = \prod_{j=1}^d f_j(x_j)$ provided that for all $j=1,\ldots,d$, $\EE \left[\prod_{k\neq j} f_k(X_k) \ind_{\left\lbrace X_k \in [a_k,b_k] \right\rbrace} \right] \neq 0$,
 \item  sums of product functions $f(\bx) =\displaystyle  \sum_{I\in \cI} \prod_{j\in I} f^I_j(x_j)$ with $\cI$ a partition of $\{1\/, \ldots \/, d\}$ provided that for all $j=1,\ldots,d$, for $ \cI \ni I \ni j$, \\
 $\EE \left[ \displaystyle\prod_{k\in I\/, \ k\neq j}f^I_k(X_k) \ind_{\left\lbrace X_k\in [a_k,b_k] \right\rbrace} \right] \neq 0$,
 \item sums of product functions $f(\bx) =\displaystyle \sum_{I\in \cI} \prod_{j\in I} f_j^I(x_j)$ with $I$ a family of subsets of  $\{1\/, \ldots \/, d\}$ provided that the $f_j^I$ are either all  increasing or all decreasing and they are either all positive or all negative,
 \item other functions, such as - for example -  in dimension $2$ with $(X_1\/,X_2) \leadsto \cU \left([0,1]^2 \right)$ - $f(x,y) = \ln (x+y)$. 
\end{enumerate}
\end{remark}
A particular case of item 4) above is the set of linear combination of Gaussian radial basis functions on $[0\/,1]^d$, with positive weights:
$$\mathcal{G} = \left\{ \sum_{i=1}^p a_i \exp[\sum_{j=1}^d (x_j-\mu_j)^2 \sigma_j^2] \/, \ a_i\geq 0\/, \ \sigma_j \geq 0\/, \ \mu_j \in \RR\right\}\/.$$
It is known that the class $\mathcal{G}$ is dense in the set of non-negative continuous functions on $[0\/,1]^d$ (see \cite{radial_functions} and also \cite{klusowki} where the class $\mathcal{G}$ is also considered). 
\begin{hyp}\label{hyp:var}\ \\
\vspace*{-0.5cm}
\begin{itemize}
 \item $Y=m(\bX)+\eps$;
 \item $\bX=(X_1\/,\ldots\/, X_d)$ is a \verob{continuous} random vector with independent coordinates;
 \item $\eps$ and $\bX$ are independent, $\eps$ is a continuous, centered random variable with increasing distribution function, $\eps$ has light tails  i.e. there exists $0<\theta<1$ such that  for any \verob{$D>0$, $\PP \left( |\eps| > D \right) \leq C \theta^D$}. 
 \item $\bX$ takes its values in $\cX$ which is assumed to be a compact hyper-rectangle of $\RR^d$: $\cX = \displaystyle \prod_{i=1}^d [u_i\/,v_i]$, $-\infty < u_i \leq v_i <\infty$;
 \item for any $y$, $\bx\mapsto F(y|\bx)$ is continuous and increasing.
\end{itemize}
\end{hyp}
\begin{theo}\label{theo:4:simple}
 Let $Y$ satisfy Assumption \ref{hyp:var}, with $m$ belonging to the $\spadesuit$-class, assume that for fixed $\beta >\frac52$, $C>0$, each constructed tree  is the highest such that $C \sqrt{n}(\ln n)^\beta\leq N_n^b\left(\Theta_\ell,\cD_n \right)$, let $k=O(n^\alpha)$, with $\alpha >0$. Then
 \[ 
\forall \bx \in \cX,  \quad \sup_{y \in \RR} \left| F_{k,n}^b \left( \left. y \right| \bX=\bx \right) -  F \left( \left. y \right| \bX=\bx \right) \right| \overset{a.s.}{\underset{n \rightarrow \infty}{\longrightarrow}} 0 \/,
\]
and
\[ 
\forall \bx \in \cX,  \quad \sup_{y \in \RR} \left| F_{k,n}^o \left( \left. y \right| \bX=\bx \right) -  F \left( \left. y \right| \bX=\bx \right) \right| \overset{a.s.}{\underset{n \rightarrow \infty}{\longrightarrow}} 0 \ .
\]
\end{theo}
\begin{remark}\label{rk:N_n}
Let us mention that the number $N_n^b\left(\Theta_\ell,\cD_n \right)$ of elements from $\cD_n^\star \left( \Theta_\ell \right)$ in each leaf may be controlled during the construction of the forest. Indeed, in most of the implementations, it is a parameter that may be chosen by the user. \\
The algorithm used to get the highest tree such that $C \sqrt{n}(\ln n)^\beta\leq N_n^b\left(\Theta_\ell,\cD_n \right)$ is as follows: choose $(j\/,z)$, with $j\in {\mathcal M}_{\mbox{try}}$ maximizing $L_A^n(j\/,z)$ while the number of elements from $\cD_n^\star \left( \Theta_\ell \right)$ in $A\cap \{[a_j\/, z]\}$ and $A\cap \{]z\/,b_j]\}$ is greater than $C \sqrt{n}(\ln n)^\beta$. \\
 In section \ref{sec:variation}, we shall see that with this procedure the height of the trees goes to infinity with $n$ a.s. (we do not have and do not need further control on the height). 
\end{remark}

It will be clear from the discussion of Section \ref{sec:variation} that Theorem \ref{theo:4:simple} is a direct consequence of Theorems \ref{theo:4:num_1} and \ref{theo:4:num_2} that we now state.  Remark that Theorems \ref{theo:4:num_1} and \ref{theo:4:num_2}  could apply in other contexts than Theorem \ref{theo:4:simple}, possibly with some dependencies, the main point being to satisfy Assumption \ref{hyp:4:num_1}. For example, for some adhoc tree constructions, Assumption \ref{hyp:4:num_1} is satisfied without independence assumption, see Remark \ref{rk:var}. In \cite{klusowki}, some hypothesis weaker than independence between the $X_i$'s are given but  the specific examples provided there are also with independent $X_i$'s.  The following assumptions allow to prove} the \verob{a.s.} consistency of our estimators in a general framework.

\begin{hyp}\label{hyp:4:num_1}
For all $\ell \in \llbracket 1,k \rrbracket$, we assume that the variation of the conditional cumulative distribution function within any cell goes to 0:
\[
\forall \bx \in \cX, \forall y \in \RR, \quad \sup_{\bz \in A_n \left( \bx;\Theta_\ell, \cD_n \right)} \left| F \left( \left. y \right| \bz \right) - F \left( \left. y \right| \bx \right) \right| \overset{a.s.}{\underset{n \rightarrow \infty}{\longrightarrow}} 0 \ .
\]
\end{hyp}
We shall discuss further on \Cref{hyp:4:num_1}. In particular in Section \ref{sec:variation}, we show that this \veroc{assumption is verified if Assumption \ref{hyp:var} is satisfied with $m$ in the $\spadesuit$-class.}  
\veroc{\begin{remark}\label{rk:var}
 Let us remark that \Cref{hyp:4:num_1} is satisfied, for example, provided that the diameter of each tree cell goes to zero and for all $y$, $F(y | \cdot)$ is continuous. This is satisfied, for example if the splitting rules in the construction of the trees imply that each direction $j=1\/, \ldots \/, d$ is chosen with positive probability at each split and a positive proportion of the sub-sample goes in each child node (see \cite{meinshausen2006quantile} Assumption 3 and Lemma 2). Imposing that each direction is chosen with a positive probability at each split seems unrealistic since even non informative variables will be chosen. This is why we propose an approach where no additional assumptions on the splitting rules are needed, see Remark \ref{rk:N_n} and Section \ref{sec:variation}.  
\end{remark}}

\begin{hyp}\label{hyp:4:num_2}
We shall make the following assumptions on $k$ (number of trees) and $N_n^b \left( \bx;\Theta,\cD_n \right)$ (number of bootstrap observations in a leaf node):
\begin{enumerate}
\item $k = \cO \left( n^\alpha \right), \textnormal { with } \alpha > 0$.
\item $\forall \bx \in \cX, \quad N_n^b \left( \bx;\Theta,\cD_n \right) = \Omega \left( \sqrt{n} \left( \ln \left( n \right) \right)^\beta \right),$ \textnormal{with} \footnote{$f \left( n \right) = \Omega \left( g \left( n \right) \right) \iff \exists k>0, \exists n_0 >0 \ \vert \  \forall n \geqslant n_0 \quad \lvert f \left( n \right) \rvert \geqslant k \cdot \lvert g \left( n \right) \rvert $} $\beta > 1$, a.s.
\item[] \textnormal{or}
\item $\forall \bx \in \cX, \quad \EE \left[ N_n^b \left( \bx;\Theta,\cD_n \right) \right] = \Omega \left( \sqrt{n} \left( \ln \left( n \right) \right)^\beta \right), \textnormal { with } \beta > 1, \textnormal{ and}$ \\
$\forall \bx \in \cX, \quad \CV \left( N_n^b \left( \bx;\Theta,\cD_n \right) \right) = \cO \left( \dfrac{1}{ \kev{n^{(\alpha + 1)/2}} \left( \ln \left( n \right) \right)^{\gamma / 2}} \right),$ \textnormal {with} \footnote{$\CV \left( X \right) = \sigma_{X} / \EE \left[ X \right]$} $\gamma> 1$.
\end{enumerate}
\end{hyp}
\begin{remark}\label{remark:4:ae_hp}
In order to prove our main consistency result, either \Cref{hyp:4:num_2} item 2. or item 3. is needed. Item 2. may seem much stronger than item 3. but it has to be noted that the number of bootstrap observations in a tree leaf is a construction parameter of the forest, so that it can be controlled. Using item 2. simplifies the proof but item 3. is sufficient.    
\end{remark}

\begin{hyp}\label{hyp:4:num_3}
For every $\bx \in \cX$, the conditional cumulative distribution function \veroc{$y \mapsto F \left( \left. y \right| \bX=\bx \right)$ is continuous and increasing.}
\end{hyp}

The two theorems below give the uniform a.s. consistency of our two estimators.
\begin{theo}\label{theo:4:num_1}
Consider a random forest which satisfies \Cref{hyp:4:num_1,hyp:4:num_2,hyp:4:num_3}. Then,
\[ 
\forall \bx \in \cX,  \quad \sup_{y \in \RR} \left| F_{k,n}^b \left( \left. y \right| \bX=\bx \right) -  F \left( \left. y \right| \bX=\bx \right) \right| \overset{a.s.}{\underset{n \rightarrow \infty}{\longrightarrow}} 0 \ .
\]
\end{theo}

\begin{theo}\label{theo:4:num_2}
Consider a random forest which satisfies \Cref{hyp:4:num_1,hyp:4:num_2,hyp:4:num_3}. Then,
\[ 
\forall \bx \in \cX,  \quad \sup_{y \in \RR} \left| F_{k,n}^o \left( \left. y \right| \bX=\bx \right) -  F \left( \left. y \right| \bX=\bx \right) \right| \overset{a.s.}{\underset{n \rightarrow \infty}{\longrightarrow}} 0 \ .
\]
\end{theo}
\begin{remark}\label{remark:4:cv_quantile}
Using standard arguments, the consistency of quantile estimates stems from \Cref{hyp:4:num_3} as well as the uniform convergence of the conditional distribution function estimators obtained above.
\end{remark}


\veroc{Let us mention that \Cref{hyp:4:num_2} allows}  to control the estimation error of our estimators and expresses that cells should contain a sufficiently large number of points so that averaging among the observations is effective.

Finally, \Cref{hyp:4:num_3} is \verob{used} to get uniform convergence of the estimators.

\veroc{Section \ref{sec:variation} is devoted to comments on Assumption \ref{hyp:4:num_1} and to prove that Assumption \ref{hyp:var} together with a control on the leave size is sufficient to have Assumption \ref{hyp:4:num_1} for functions of the $\spadesuit$-class. }


\section{On variation in Random Forest leaves}\label{sec:variation}
\veroc{
This section is devoted to some further analysis on the variation of $F(y|\cdot)$ on leaves, as formalized in Assumption \ref{hyp:4:num_1}. 
\subsection{Some comments on Assumption \ref{hyp:4:num_1}}
\Cref{hyp:4:num_1} ensures a control on the approximation error of the estimators. It is drawn \verob{from  Proposition 2 of \cite{scornet2015consistency} which} shows the consistency of Breiman's random forest estimate in an additive regression framework. Their Proposition 2 allows to manage the approximation error of the estimator by showing that the variation of the regression function $m$ within a cell of a random empirical tree is small provided $n$ is large enough. This result is based on \verob{ Lemma 1 of \cite{scornet2015consistency}} which states that the variation of the regression function $m$ within a cell of a random theoretical tree goes to zero for an additive regression model. A random theoretical tree is grown \verob{following the same rules as a} random empirical tree, except that the theoretical equivalent of the empirical CART-split criterion (\ref{def:L_n}) \verob{on a} node A below is used to choose the best split
\begin{align}\label{eq:theo_split}
L_A^\star (i,z) = & \var \left( \left. Y \right| \bX \in A \right) \\
& - \PP \left( \left. X_i < z  \right| \bX \in A \right) \var \left( \left. Y \right| X_i < z , \bX \in A \right) \nonumber \\
& - \PP \left( \left. X_i \geqslant z  \right| \bX \in A \right) \var \left( \left. Y \right| X_i \geqslant z , \bX \in A \right) \ . \nonumber
\end{align}
Hence, a theoretical tree is obtained thanks to the best consecutive cuts $(i^\star, z^\star)$, among $i \in {\mathcal M}_{\mbox{try}}$, $z\in A^i$ optimizing the previous criterion $L^\star (\cdot, \cdot)$. Remark that ${\mathcal M}_{\mbox{try}}$ is selected with respect to $\Theta^2$, as for the empirical forest. 

General results on standard partitioning estimators whose construction is independent of the label in the training set (see Chapter 4 in \cite{gyorfi2006distribution} or Chapter 6 in \cite{devroye2013probabilistic}) state that a \vero{sufficient} condition to prove the consistency is that the diameter of the cells tend to zero as $n \rightarrow \infty$. Instead of such a geometrical assumption, Proposition 2 in \cite{scornet2015consistency} ensures that the variation of $m$ inside a node is small thanks to their Lemma 1. But the cornerstone of the Lemma 1 is the Technical Lemma 1 of \cite{scornet2015supplementary} recalled below for completeness.
\begin{lem*}[Technical]\label{lem:4:technical_lemma}
Assume that: 
\begin{itemize}
\item $Y = m \left( \bX \right) + \varepsilon$ with $m \left( \bX \right) = \sum\limits_{i=1}^d m_i \left( X_i \right), \bX \sim \cU \left( \left[ 0,1 \right]^d \right)$ and $\varepsilon \sim \cN(0, \sigma^2)$,
\item $L_A^\star (i,z)=0 \quad \verob{\forall i = 1\/, \ldots\/,d} , \forall z \in \left[ a_i, b_i \right] \  \left( 0 \leqslant a_i < b_i \leqslant 1 \right)$,
\end{itemize}
then the regression function $m$ is constant on $A$.
\end{lem*}
This lemma states that if the theoretical split criterion is zero for all cuts in a node, then the regression function $m$ is constant on this node, i.e. the variation of \verob{$m$ on}  the cell is zero. But, \vero{examples for which $L_A^\star (i,z)=0 \quad \forall i, \forall z \in \left[ a_i, b_i \right]$ and the regression function is not constant \kev{can be} easily constructed.} 
\vero{Let us consider a two-dimensional example, let $A = A_1 \times A_2 = \left[ a_1, b_1 \right] \times \left[ a_2, b_2 \right]$ and suppose that the response $Y$ is}
 
\begin{equation}\label{eq:counter}
 Y = X_1 X_2 + c_1 X_1 + c_2 X_2 + \varepsilon \  =: m(X_1\/,X_2) + \varepsilon\/,
 \end{equation}
 
 with

\begin{itemize}
\item $\bX = \left( X_1, X_2 \right)$ independent random inputs,
\item $c_1 = -\dfrac{\EE \left[ X_2 \ind_{\left\lbrace X_2 \in A_2 \right\rbrace} \right]}{\PP \left( X_2 \in A_2 \right)}$ and $c_2 = -\dfrac{\EE \left[ X_1 \ind_{\left\lbrace X_1 \in A_1 \right\rbrace} \right]}{\PP \left( X_1 \in A_1 \right)}$,
\item and $\varepsilon$ a centered noise independent of $\bX$.
\end{itemize}
It can be shown for this model that within the node $A$, $L^\star \equiv 0$ for all $i \in \left\lbrace 1,2 \right\rbrace$,  for all $z \in \left[ a_i, b_i \right]$ and yet the regression function $m$ is not constant. \\
\ \\
Accordingly, the technical lemma above is well-designed for an additive regression framework. But this context is far from reality for many concrete examples.
\subsection{Sufficient conditions to control the variation on leaves}}
Our purpose is \veroc{now} to prove that \veroc{Assumption \ref{hyp:var} with a control on the size of the leaves and a regression function in the $\spadesuit$-class is sufficient to insure \Cref{hyp:4:num_1} which controls  the variation of $F(y|\cdot)$ on leaves.} An interesting sub-product of this study is Proposition \ref{prop:close} which shows the asymptotic proximity of empirical and theoretical trees. \\
%
In what follows, \veroc{$C$ denotes any positive constant so that we may write $C+C=C$, $uC=C$, with $u>0$, ....} 

We are interested in $F(y|\cdot)$, remark that $F(y|\cdot)$ is constant on $A$ if and only if $m(\cdot)$ is constant on $A$. Indeed,
\[
F(y|\bx) = \PP(Y\leq y | \bX=\bx) = \PP(\eps\leq y-m(\bx)) = F_\eps(y-m(\bx))\ .
\]
We shall use \verob{ theoretical trees} construction as in \cite{scornet2015consistency},  the theoretical CART-split is done by maximizing $L_A^\star (i,z)$ defined in \eqref{eq:theo_split}. 
If $L_A^\star$ admits several maxima, then one is chosen uniformly at random. \\
Let us denote by  $A_n\left( \Theta_\ell,\cD_n \right)$ any leaf \kev{in} an empirical tree, and $N_n^b \left(\Theta_\ell,\cD_n \right)$ be the number of elements of $\cD_n^\star \left( \Theta_\ell \right)$ that fall into $A_n \left(\Theta_\ell,\cD_n \right)$. We  have the following result.
\veroc{
\begin{theo}\label{theo:variation}
Let $Y$ satisfy \Cref{hyp:var}, with $m$ belonging to the $\spadesuit$-class, let $\beta >\frac52$, $C >0$, let the constructed trees be the highest such that $C \sqrt{n}(\ln n)^\beta\leq N_n^b\left(\Theta_\ell,\cD_n \right) $, then \Cref{hyp:4:num_1} is verified. 
\end{theo}
}
Note that for additive $m$ it is proven in  \cite{scornet2015consistency} that the variation of $m$ on leaves goes to zero in probability, with an assumption on the height of trees. We extend this result to the $\spadesuit$-class and obtain an a.s. convergence. \verob{The control on the height of trees is replaced by a control on the size of leaves.} The strategy of proof is partially inspired by \cite{scornet2015supplementary} and \cite{scornet2015consistency}. 
%
%
The following lemma proves \veroc{that the variation of $ F(y|\cdot)$ on leaves of theoretical trees goes to $0$ a.s. as the height of the trees goes to infinity. }
\begin{lem}\label{lem:variation}
 Assume that \veroc{\Cref{hyp:var} is satisfied with  the regression function $m$ in} the $\spadesuit$-class, let $S^\infty =(s_j\/, j=1\/,\ldots)$ with $s_j\in \{L\/,R\}$, it describes an infinite path in a \veroc{binary} tree, let $S^h=(s_j\/, j=1\/,\ldots\/,h)$, it describes a path in a \veroc{binary} tree of height $h$. Let $A_h(S^h\/,\Theta)$ be the corresponding leaf in a theoretical tree. Then the variation of $F(y|\cdot)$ on $A_h(S^h\/,\Theta)$ goes to $0$ a.s. as $h$ goes to infinity. 
\end{lem}
\begin{proof}
 Since $\bx\mapsto F(y|\bx)$ is assumed to be continuous, then the result holds if the diameter of $A_h(S^h\/,\Theta)$ goes to zero. \veroc{Let} $A_\infty(S^\infty\/,\Theta) = \displaystyle\bigcap_{h\geq 1} A_h(S^h\/,\Theta) $, it is a decreasing intersection of rectangles, if its diameter is non zero then it writes: for a non empty subset $ J\subset \{1\/, \ldots \/,d\}$, let $J^c$ be its complementary, $\left(x_j\/, j\in J^c\/, a_j\/, j\in J \right)$; $\left(x_j \/, j\in J^c\/, b_j\/, j\in J\right) \in \cX$
$$A_\infty(S^\infty \/,\Theta)= \{x_j \/, j\in J^c\} \times \prod_{j\in J} [a_j\/,b_j]=\{x_j \/, j\in J^c\} \times A^J\/, \ \mbox{with} \ a_j<b_j.$$
 Remark that for any rectangle \verob{$A = \prod_{k=1}^d [\alpha_k\/, \beta_k]$},
 \begin{eqnarray*}
  L_A^\star(i\/,z) &=&  \PP \left( \left. X_i < z  \right| \bX \in A \right) \left( \EE \left[ \left. Y \right| X_i < z , \bX \in A \right] - \EE \left[ \left. Y \right| \bX \in A \right] \right)^2  \\
  &+& \PP \left( \left. X_i \geqslant z  \right| \bX \in A \right) \left( \EE \left[ \left. Y \right| X_i \geqslant z , \bX \in A \right] - \EE \left[ \left. Y \right| \bX \in A \right] \right)^2 
   \end{eqnarray*}
 \verob{can be seen as a continuous function of $\left( \alpha_k, \beta_k \right)_{k=1, \ldots, d}, j\/, z$. \\
 For $j\in J$,} $L^\star_{A_\infty \left(S^\infty, \Theta \right)}(j\/,z)$ rewrites
 \begin{eqnarray*}
  L_{A_\infty \left( S^\infty\/,\Theta \right)}^\star(j\/,z) &=&  \PP \left( \left. X_j < z  \right| \bX^J \in A^J \right) \left( \EE \left[ \left. Y \right| X_j < z , \bX^J \in A^J \right] - \EE \left[ \left. Y \right| \bX^J \in A^J \right] \right)^2  \\
  &+& \PP \left( \left. X_j \geqslant z  \right| \bX^J \in A^J \right) \left( \EE \left[ \left. Y \right| X_j \geqslant z , \bX^J \in A^J \right] - \EE \left[ \left. Y \right| \bX^J \in A^J \right] \right)^2\/.
   \end{eqnarray*}
   Also, for $j\in J^c$, $ L_{A_\infty(S^\infty)}^\star(j\/,z)   =0$. \\
 Write $A_h(S^h\/,\Theta) =\prod_{j=1}^d [a_j^h\/, b_j^h]$, then for $j \in J$, $a_j^h \rightarrow a_j$ and $b_j^h\rightarrow b_j$ as $h \rightarrow \infty$, also $L_{A_h(S^h)}^\star(j\/,z) \rightarrow  L_{A_\infty(S^\infty)}^\star(j\/,z)$.  \\
 Let $(j^\star,z^\star) \in \mbox{argmax}\ L_{A_\infty(S^\infty)}^\star(j\/,z) $, if $j^\star \in J^c$ then $ L_{A_\infty(S^\infty)}^\star(j^\star\/,z^\star)   =0$ and thus $ L_{A_\infty(S^\infty)}^\star(j\/,z)   =0$ for any $j,z$. Assume $j^\star \in J$  and consider a subsequence $(h_p)_{p\in \NN}$ such that $j^\star \in \cM_{\mbox{try}}$ at each level $h_p$, it exists a.s. because $\Theta$ gives a positive probability \verob{ to each $j \in \{1\/,\ldots\/,d\}$ of belonging to $\cM_{\mbox{try}}$} at each level. Then, 
 consider $(j^p\/,z^p) \in \mbox{argmax}\ L_{A_{h_p}(S^{h_p})}^\star(j\/,z) $ and $(j^\infty\/,z^\infty)$ any limit point of the sequence $(j^p\/,z^p)$, which means that for a subsequence $p_q \rightarrow \infty$, $(j^{p_q}\/, z^{p_q}) \rightarrow (j^\infty\/, z^\infty)$, thus for $q$ large enough, $j^{p_q}=j^\infty$. It easily follows that $ L_{A_\infty(S^\infty)}^\star(j^\star\/,z^\star)    \leq  L_{A_\infty(S^\infty)}^\star(j^\infty\/, z^\infty)$ and thus $(j^\infty\/, z^\infty)\in \mbox{argmax}\ L_{A_\infty(S^\infty)}^\star(j\/,z)$. Now, if $j^\infty \in J$, each $z^h$ is either $a_{j_h}^h$ or $b_{j_h}^h$; $z^\infty$ is either $a_{j^\infty}$ or $b_{j^\infty}$. So that either $\{X_{j^\infty} < z , \bX^J \in A^J \}= \{X^J\in A^J\}$ and  $\{X_{j^\infty} \geq z , \bX^J \in A^J \}= \emptyset$ or $\{X_{j^\infty} < z , \bX^J \in A^J \}= \emptyset$ and  $\{X_{j^\infty} \geq z , \bX^J \in A^J \}= \{X^J\in A^J\}$, thus
 \begin{eqnarray*}
  L_{A_\infty(S^\infty)}^\star(j^\infty\/,z^\infty)  &=&  \PP \left( \left. X_{j^\infty} < z  \right| \bX^J \in A^J \right) \left( \EE \left[ \left. Y \right| X_{j^\infty} < z , \bX^J \in A^J \right] - \EE \left[ \left. Y \right| \bX^J \in A^J \right] \right)^2  \\
  &+& \PP \left( \left. X_{j^\infty} \geqslant z  \right| \bX^J \in A^J \right) \left( \EE \left[ \left. Y \right| X_{j^\infty} \geqslant z , \bX^J \in A^J \right] - \EE \left[ \left. Y \right| \bX^J \in A^J \right] \right)^2\/\\
  &=&0\/.
 \end{eqnarray*}
If $j^\infty \in J^c$, then \verob{$  L_{A_\infty(S^\infty)}^\star(j^\infty\/,z^\infty)=0$}. Finally, since $ L_{A_\infty(S^\infty)}^\star(j\/,z) \leq  L_{A_\infty(S^\infty)}^\star(j^\infty\/,z^\infty)=0$, we conclude that $L_{A_\infty(S^\infty)}^\star (j,z)=0$ for all $(j,z)$. This is equivalent to: for all $i=1\/,\ldots \/,d$, $z\in [a_i\/,b_i]$,
\begin{eqnarray*}
\lefteqn{\EE [  Y | X_i \leq z , \bX \in A_\infty ] - \EE [  Y | \bX \in A_\infty ]= 0 \Leftrightarrow} \\
&&\PP( \bX \in A_\infty) \EE \left[ Y \ind_{\left\lbrace X_i \leq z, \ \bX \in A_\infty \right\rbrace} \right] = \PP( X_i \leq z , \bX \in A_\infty) \EE \left[ Y\ind_{\left\lbrace \bX\in A_\infty \right\rbrace} \right].
\end{eqnarray*}
By derivating with respect to $z$, we may see that it is equivalent to $z\mapsto \EE \left[ m \left( z, \bX^{-i} \right) \ind_{\left\lbrace \bX^{-i}\in A_\infty^{-i} \right\rbrace} \right]$ is constant for all $i=1,\ldots,d$. Since we assumed that $m$ belongs to the $\spadesuit$-class, either $m$ is constant on $ A_\infty(S^\infty)$ or the diameter of $A_\infty(S^\infty)$ is zero. In both cases, the variation of $F(y|\cdot)$ on $A_h(S^h\/,\Theta)$ goes to $0$ as $h$ goes to infinity. 
\end{proof}
 \veroc{We now show that empirical and theoretical trees are close.}  The main tool for this step is the following result which has an intrinsic interest since it proves that the theoretical and empirical cost functions used for the tree constructions are uniformly close, provided that the leave sizes are not too small. Recall that $L_A^n(j\/,z)$ is the empirical cost on the rectangle $A$ on direction $j$ and at level $z$. 
\begin{prop}\label{prop:Lstar}
Let \Cref{hyp:var} be satisfied. Let $\beta >\frac52$, let $A$ be a rectangle in $\cX$, we shall say that $(A\/,j\/,z) \in \cA^n$ if  the numbers $N_{A_L}^b$, $N_{A_R}^b$ of elements of $\cD_n^\star$ belonging to $A_L:= A \cap \{x_j \leq z\}$ and $A_R := A\cap \{x_j>z\}$ \verob{are greater} than $C \sqrt{n}(\ln n)^\beta$. We have
\[
\sup_{(A\/,j\/,z)\in \cA^n} | L_A^\star(j\/,z) - L_A^n(j\/,z)| \stackrel{\mbox{a.s.}}{\longrightarrow} 0 \ \mbox{as} \ n \rightarrow \infty\/.
\]
\end{prop}
\begin{proof}
Let $(A\/,j\/,z)\in \cA^n$ be fixed. Rewrite the difference $| L_A^\star(j\/,z) - L_A^n(j\/,z)| =  |T_L+T_R|$ with
\begin{eqnarray*}
T_L &=& T_{L,1} + T_{L,2} \\
 	&=:& \dfrac{N_{A_L}^b}{N_A^b} \left( \left(\overline{Y_A}-\overline{Y_{A_L}} \right)^2 - \left( \EE \left[Y|\bX\in A\right] - \EE \left[ Y|\bX \in A_L \right] \right)^2\right) \\
	&& + \left( \EE \left[Y|\bX\in A \right] - \EE \left[ Y|\bX\in A_L \right] \right)^2 \left( \frac{N_{A_L}^b}{N_A^b}- \PP \left( \bX\in A_L|\bX\in A \right)\right) \ .
\end{eqnarray*}
 The term $T_R$ is defined in the same way by using $A_R$ instead of $A_L$. In order to prove the proposition, we shall prove that $\sup_{A\/,j\/,z} T_{L\/,1}$ and $\sup_{A\/,j\/,z} T_{L\/,2}$ go to $0$ a.s. The same holds in the same way for $T_{R\/,1}$ and $T_{R\/,2}$. Using Vapnik-Chervonenkis theory on rectangles in $\cA$ (see Lemma \ref{lem:vc_app} in Appendix \ref{sec:proof2}) we have:
 \begin{equation}\label{eq:vc1}
 \PP\left(\sup_{B\in\cA}\left|\frac{N_B^b}n - \PP(\bX\in B)\right|>\kappa\right)\leq 16 \kev{(n+1)}^{2d}e^{-n\kappa^2/\kev{128}}\/.
 \end{equation}
$T_{L\/,2}$ decomposes into:
$$|T_{L\/,2}| \leq (\EE(Y|\bX\in A)-\EE(Y|\bX\in A_L))^2 \times \left[\frac{n}{N_A^b}\left|\frac{N_{A_L}^b}n- \PP(\bX\in A_L)\right| + \left|\frac1{\PP(\bX\in A)}-\frac{n}{N_A^b}\right|\PP(\bX \in A_L) \right]\/.$$
Remark that for $B \in \cA$, if
\begin{equation}\label{g_proba}
 \left| \frac{N_B^b}n - \PP(\bX \in B) \right| \leq \frac{C}2 \frac{(\ln n)^\beta}{\sqrt{n}}
\end{equation}
and $N_B^b \geq C \sqrt{n}(\ln n)^\beta$, then $\PP(\bX \in B) \geq   \frac{C}2 \frac{(\ln n)^\beta}{\sqrt{n}}$. So that, for $(A\/,j\/,z) \in \cA^n$, we have, provided that (\ref{g_proba}) holds for $A$,
\veroc{
\begin{eqnarray*}
 \EE(Y|\bX \in A) &=& \frac1{\PP(\bX \in A)} \left[\EE(Y\ind_{\{\bX\in A\}}\ind_{\{Y\leq D\}}) + \EE(Y\ind_{\{\bX\in A\}}\ind_{\{Y> D\}})\right]\\
 &\leq& D + \EE(Y^p)^\frac1p \frac{\PP(Y>D)^\frac1q}{\PP(X\in A)^{1-\frac1r}}\\
 \lefteqn{\mbox{the second term is obtained using Hölder inequality}}\\
 &\leq & (\ln n)^\gamma + C e^{(\ln n)^\gamma \frac{\ln \theta}q} \frac{n^{\frac12 (1-\frac1r)}}{(\ln n)^{\beta(1-\frac1r)}} \\
 & \leq & C(\ln n)^\gamma\/,
\end{eqnarray*}
}
by taking $D = (\ln n )^\gamma$, $\gamma >1$, $p\/,q\/,r>0$ with $\frac1p+\frac1q+\frac1r = 1$.   Now, Equation \eqref{eq:vc1} gives that \eqref{g_proba} is satisfied for $A$ and $A_L$ with probability greater than $1- \kev{16(n+1)}^{2d}e^{-C\frac{(\ln n)^{2\beta}}{\kev{512}}}$ and the condition on $N_{A_L}$ for $(A\/,j\/,z) \in \cA^n$ gives that
\[
\PP(T_{L\/,2} >\kappa)\leq \kev{16(n+1)}^{2d}e^{-C\frac{(\ln n)^{2\beta}}{\kev{512}}} + C \kev{(n+1)}^{2d} e^{-C \kappa^2 \frac{(\ln n)^{2\beta}(\ln n)^{4\gamma}}{128}} \ .
\]
Then, Borel-Cantelli Lemma gives that $\sup_{(A\/,j\/,z)\in \cA^n} T_{L\/,2}$ goes to $0$ a.s. provided that $2\beta -4\gamma >1$. \\
The term $T_{L\/,1}$ is treated in the same way, instead of using (\ref{eq:vc1}), we use Lemma \ref{lem:giorfy} in Appendix \ref{sec:proof2} to get that for any \verob{$D>0$}, for any $\kappa>0$, $A\in \cA$ and $\frac1p+\frac1q=1$
\verob{\begin{eqnarray}\label{eq:gyorfi}
 \lefteqn{\PP\left(\left|\frac1n \sum_{i=1}^nY_i^\star \ind_{\{\bX_i^\star\in A\}} - \EE(Y\ind_{\{\bX \in A\}})\right| >\kappa \right)} \nonumber \\
 &\leq& 6\left(\frac{24eD}\kappa \ln \left(\frac{48eD}\kappa\right)\right)^{2d} e^{-n\kappa^2/(128 D^2)} + C\frac{\EE(Y^p)^\frac1p \PP(Y>D)^\frac1q}{\kappa}\/,
\end{eqnarray}
where the $(\bX_i^\star\/,Y_i^\star)$'s  \veroc{form a bootstrap sample} from $(\bX_1\/,Y_1)\/, \ldots \/, (\bX_n\/,Y_n)$ and we take $D=(\ln n)^\gamma$ as before. }
\end{proof}
The last stone for the proof of Theorem \ref{theo:variation} is to prove that \verob{at} each level $h$, each node of the empirical tree $A^n(S)$ is close to a level $h$ node of a theoretical tree. 
\begin{prop}\label{prop:close}
Let Assumption \ref{hyp:var} be satisfied. Assume that for $\beta >\frac52$, $ N_n^b\left(\Theta_\ell,\cD_n \right) \geq C \sqrt{n}(\ln n)^\beta$.
 For $h\in \NN$, let $S\in \{L\/,R\}^h$ describe a path of length $h$ in a binary tree, let $A^n(S)$  and $A(S)$ be  corresponding nodes in  empirical and theoretical trees. Denote
 $$A(S) = \prod_{j=1}^d [a_j\/, b_j] \ \mbox{and} \ A^n(S) =  \prod_{j=1}^d [a_j^n\/, b_j^n]\/.$$
 \veroc{Denote $\mathcal{T}_h$ the set of theoretical trees of height $h$, then
 \begin{equation}\label{eq:close}\inf_{\mathcal{T}_h} \max_{j=1\/, \ldots \/, d} \max \left(|a_j-a_j^n|\/, |b_j-b_j^n|\right) \longrightarrow 0 \  \mbox{a.s.} \ \mbox{as} \ n \rightarrow \infty
 \end{equation}
 }
\end{prop}
\veroc{\begin{remark}
Proposition \ref{prop:close} implies that for any $h\in \NN$ fixed, for any empirical tree with node sizes greater than $C\sqrt{n}(\ln n)^\beta$, we may find a theoretical one as close as we wish, until height $h$. This is an a.s. and more precise version of Lemma 3 in \cite{scornet2015consistency}.
\end{remark}
}
\begin{proof}
 We proceed by induction on $h$. If $h=0$ then $A(S) = A^n(S) = \cX$ and the assertion holds. Let us assume that the result holds for $h \in \NN$, let $S'\in  \{L\/,R\}^{h+1}$ describe a path of length $h+1$ in a binary tree, \verob{$S' = (S \/, u)$} with $ S\in\{L\/,R\}^h$ and $u\in \{L\/,R\}$. \verob{ Let $(j^n\/, z^n) \in \mbox{argmax}\ L^n_{A^n(S)}(j\/,z)$, then
 $$ A^n(S') = A^n(S) \cap \{x_{j^n} \leq z^n\} \ \mbox{or} \  A^n(S') = A^n(S) \cap \{x_{j^n} > z^n\}\/. $$
 Let us denote $\displaystyle A^n(S) = \prod_{j=1}^d[a_j^n\/, b_j^n]$. }
  \verob{Our} hypothesis on the construction of the empirical tree implies that \verob{for all $n\in \NN$,}  $(A^n(S)\/,j^n\/,z^n) \in \cA^n$ so that by Proposition \ref{prop:Lstar}:
\begin{equation}\label{eq:proche}
|L^n_{A^n(S)} (j^n\/,z^n) - L^\star_{A^n(S)}(j^n\/, z^n)| \longrightarrow 0 \ \mbox{a.s. as} \ n \rightarrow \infty \ .
\end{equation}

 The induction hypothesis implies that, \veroc{taking if necessary a subsequence, we may find $A(S) = \displaystyle\prod_{j=1}^d [a_j\/,b_j]$  a node in a theoretical tree at height $h$ such that $|a_j^n-a_j|\rightarrow 0$ and $| b_j^n-b_j| \rightarrow 0$ a.s.} Let $(j^\infty\/, z^\infty)$ be any limit point of $(j^n\/, z^n)$, the continuity of $L^\star$ as a function of $(a_j\/, b_j\/, j\/, z)$ and (\ref{eq:proche}), lead to
%
%
$|L^{n}_{A^n(S)}(j^n\/, z^n) - L^\star_{A(S)}(j^\infty\/, z^\infty)| \longrightarrow 0$, a.s. as $n$ goes to infinity, taking if needed a subsequence. \\
\veroc{Let us consider any $(j\/,z)$,  if $ (A^n(S)\/, j\/, z) \not\in \cA^n$ for $n$ large enough, then  $L^\star_{A(S)}(j\/,z) = 0$ because, in that case, either:\\
$\PP( X_j < z \/, \bX \in A(S))=0$ and $\PP( X_j \geq z,  \bX\in A(S))=\PP(\bX \in A(S))$ or\\
$\PP( X_j \geq z \/,  \bX\in A(S))=0$ and $\PP( X_j < z \/, \bX \in A(S))=\PP(\bX \in A(S))$ (using (\ref{eq:vc1})).}\\ Otherwise, $ (A^n\/, j\/, z) \in \cA^n$ for infinitely many $n$ and using again  Proposition \ref{prop:Lstar}, we have  $L^\star_{A(S)}(j^\infty\/, z^\infty) \geq L^\star_{A(S)}(j\/,z)$.  We conclude that $(j^\infty\/, z^\infty) \in \mbox{argmax}\ L^\star_{A(S)}(j\/,z)$ and thus 
$$ A(S) \cap \{x^{j^\infty} \leq z^\infty\} \ \mbox{and} \ A(S) \cap \{x^{j^\infty} > z^\infty\}$$
are level $h+1$ nodes of a theoretical tree. In other words, $z^n$ and $z^\infty$ are respectively new end points of $A^n(S')$, $A(S')$ and since \veroc{$ z^\infty$ is a limit point of $(z^n)_{n\in\NN}$}, this concludes the induction step and the proof. 
\end{proof}
\veroc{
\begin{proof}[End of proof of Theorem \ref{theo:variation}]
In order to conclude the proof of Theorem \ref{theo:variation}, we need to prove that  the variation of $F(y|\cdot)$ on any leaf in an empirical tree goes to $0$ a.s.  \\
Consider $h\in \NN$, a leaf is described by a finite sequence $S\in \{L\/, R\}^h$ of some $S^\infty \in \{L\/,R\}^\NN$. If $A^n(S)$ is a node of the empirical tree, it means that  the number of elements from $\cD_n^\star$ in $A^n(S)$  is greater than $C \sqrt{n}(\ln n)^\beta$. \\
Let $(j^n\/, z^n) \in \mbox{argmax}\ L^n_{A^n(S)}(j\/,z)$. On one side, if $(A^n(S)\/,j^n\/, z^n ) \not\in \cA^n$ for infinitely many $n$, let $A(S)$ be a node in a theoretical tree, given by Proposition \ref{prop:close}. As in the proof of Proposition \ref{prop:close},  we have that $L^\star_{A(S)}(j\/,z)=0$ for all $(j\/,z)$, since $m$ belongs to the $\spadesuit$-class, this implies that the variation of $F(y|\cdot)$ on $A(S)$ is zero. Using again Proposition \ref{prop:close}, we conclude that for any $\kappa>0$, for $n$ large enough, the variation of $F(y|\cdot)$ on $A^n(S)$ is less than $\kappa$. In this case, $A^n(S)$ is a leaf of the empirical tree (since the condition on the number of elements in the child nodes is not satisfied). \\
Otherwise,  $(A^n(S)\/,j^n\/, z^n ) \in \cA^n$ for $n$ large enough, which means that the construction of the empirical tree continues. \\
By induction on $h$, we get that if $A^n(S)$ is a leaf in an empirical tree, with $|S|=h_n$ then either the variation of $F(y|\cdot)$ on $A^n(S)$ goes to $0$ or $h_n$ goes to infinity. Consider the second case.
%
%
Fix $\kappa >0$ and using Lemma \ref{lem:variation}, choose $h$ so that for $|S|=h$ the variation of $F(y|\cdot)$ on a node $A(S)$  in an empirical tree  is less than $\kappa$. Choose $\eta$ so that if $|z-z'|<\eta$ then $|F(y|z)-F(y|z')|<\kappa$. Using Proposition \ref{prop:close}, choose $n$ such that for any $z\in A^n(S)$ there is $z' \in A(S)$ with $|z-z'|<\eta$. Then, the variation of $F(y|\cdot)$ on $A_n(S)$ is less than $3\kappa$. 
\end{proof}
To conclude this section, let us say that \Cref{theo:variation} gives simple and realistic conditions under which the hypothesis of \Cref{theo:4:num_1,theo:4:num_2} are verified. 
}

The next section is devoted to the proof of  \Cref{theo:4:num_1}. The proof of \Cref{theo:4:num_2} is similar and left to the reader.

\section{Proofs of Theorem \ref{theo:4:num_1}}\label{sec:5:proofs}

The proofs of \Cref{theo:4:num_1,theo:4:num_2} are close. We \kev{only provide that of \Cref{theo:4:num_1}} \verob{below}.

The main ingredient of the proof is to use a second sample $\cD_n^{\diamond}$ in order to deal with the data-dependent aspect. Thus, we first define a dummy estimator based on two samples $\cD_n$ and $\cD^{\diamond}_n$ which will be used below. The trees are grown \vero{ using $\cD_n$}, but we consider another sample $\cD_n^{\diamond}$ (independent of $\cD_n$ and $\Theta$) which is used to define a dummy estimator
\begin{equation}\label{eq:5:CDF_estimate_diam}
F_{k,n}^{\diamond} \left( \left. y \right| \bX=\bx;\Theta_1,\ldots,\Theta_k, \cD_n^{\diamond}, \cD_n \right) = \sum_{j=1}^{n} w_{n,j}^{\diamond} \left( \bx; \Theta_1,\ldots,\Theta_k, \bX^{\diamond 1}, \ldots, \bX^{\diamond n}, \cD_n \right) \ind_{ \{Y^{\diamond j} \leqslant y\}} \ ,
\end{equation}
where the weights are
\[
w_{n,j}^{\diamond}\left( \bx; \Theta_1,\ldots,\Theta_k, \bX^{\diamond 1}, \ldots, \bX^{\diamond n}, \cD_n \right) = \dfrac{1}{k} \sum_{\ell=1}^{k} \dfrac{\ind_{\left\lbrace \bX^{\diamond j} \in A_n \left( \bx;\Theta_\ell,\cD_n \right) \right\rbrace}}{N_n^{\diamond} \left( \bx;\Theta_\ell, \bX^{\diamond 1}, \ldots, \bX^{\diamond n}, \cD_n \right)}, \ j=1,\ldots, n,
\]
with $N_n^{\diamond} \left( \bx;\Theta_\ell, \bX^{\diamond 1}, \ldots, \bX^{\diamond n}, \cD_n \right)$, the number of elements of $\cD_n^{\diamond}$ that fall into $A_n \left( \bx;\Theta_\ell,\cD_n \right)$. Throughout this section, we shall use the convention $\frac00=0$ in case 
$N_n^{\diamond} \left( \bx;\Theta_\ell, \bX^{\diamond 1}, \ldots, \bX^{\diamond n}, \cD_n \right)=0$ and thus $\ind_{\left\lbrace \bX^{\diamond j} \in A_n \left( \bx;\Theta_\ell,\cD_n \right) \right\rbrace}=0$ for $j=1\/,\ldots n$. \\
The weights $w_{n,j}^{\diamond}\left( \bx; \Theta_1,\ldots,\Theta_k, \bX^{\diamond 1}, \ldots, \bX^{\diamond n}, \cD_n \right)$ are nonnegative random variables, as function of $\Theta_1,\ldots,\Theta_k, \bX^{\diamond 1}, \ldots, \bX^{\diamond n},\cD_n$. To lighten the notation in the sequel, we will simply write $F_{k,n}^{\diamond} \left( \left. y \right| \bX=\bx \right) = \sum\limits_{j=1}^{n} w_{j}^{\diamond} \left( \bx \right) \ind_{ \{Y^{\diamond j} \leqslant y\}}$ instead of \eqref{eq:5:CDF_estimate_diam}.

Let $\bx \in \cX$ and $y \in \RR$, we have
\begin{eqnarray*}
\left| F_{k,n}^b \left( \left. y \right| \bX=\bx \right) -  F \left( \left. y \right| \bX=\bx \right) \right| & \leqslant & \left| F_{k,n}^{\diamond} \left( \left. y \right| \bX=\bx \right) -  F \left( \left. y \right| \bX=\bx \right) \right| \\
&& + \left| F_{k,n}^{\diamond} \left( \left. y \right| \bX=\bx \right) - F_{k,n}^b \left( \left. y \right| \bX=\bx \right) \right| \ .
\end{eqnarray*}
The convergence of the two right-hand terms is handled separately into the following \Cref{prop:5:cv_diam_to_true} and \Cref{lem:5:cv_diam_to_bootstrap}.
\begin{prop}\label{prop:5:cv_diam_to_true}
Consider a random forest which satisfies \Cref{hyp:4:num_1,hyp:4:num_2}. Then,
\[ 
\forall \bx \in \cX, \forall y \in \RR, \quad F_{k,n}^{\diamond} \left( \left. y \right| \bX=\bx \right) \overset{a.s.}{\underset{n \rightarrow \infty}{\longrightarrow}} F \left( \left. y \right| \bX=\bx \right) \ .
\]
\end{prop}
Hence, \Cref{prop:5:cv_diam_to_true} establishes the consistency for a random forest estimator based on a second sample $\cD_n^\diamond$ independent of $\cD_n$ and $\Theta$. \cite{wager2018estimation} proved that estimators built from honest forests are asymptotically Gaussian. Remark that in \cite{wager2018estimation}, it is also required to control the proportion of chosen observations at each split and in each direction. In our case, going through a kind of honest trees is just a theoretical tool. We go one step further with the following lemma by showing that the estimators built with \veroc{$\cD_n^*$ and $\cD_n^\diamond$} are close.
\begin{lem}\label{lem:5:cv_diam_to_bootstrap}
Consider a random forest which satisfies \Cref{hyp:4:num_2}. Then,
\[ 
\forall \bx \in \cX, \forall y \in \RR, \quad \left| F_{k,n}^{\diamond} \left( \left. y \right| \bX=\bx \right) -  F_{k,n}^b \left( \left. y \right| \bX=\bx \right) \right| \overset{a.s.}{\underset{n \rightarrow \infty}{\longrightarrow}} 0 \ .
\]
\end{lem}

Hence, according to \Cref{prop:5:cv_diam_to_true} and \Cref{lem:5:cv_diam_to_bootstrap}, we get
\begin{equation}\label{eq:5:theo1_cv_CDF_estimate_to_true_pointwise_ps}
\forall \bx \in \cX, \forall y \in \RR, \quad F_{k,n}^b \left( \left. y \right| \bX=\bx \right) \overset{a.s.}{\underset{n \rightarrow \infty}{\longrightarrow}} F \left( \left. y \right| \bX=\bx \right) \ .
\end{equation}

Now, thanks to Dini's second theorem, let us sketch how to obtain the almost sure uniform convergence relative to $y$ of the estimator.

Note that $\left\lbrace Y^{j} \leqslant y \right\rbrace = \left\lbrace U_j \leqslant F_{Y|\bX=\bx} \left( y \right) \right\rbrace$ under \Cref{hyp:4:num_3} with $U_j = F_{Y|\bX=\bx} \left( Y^{j} \right), j=1,\ldots,n$ which are i.i.d random variables. Then, \eqref{eq:5:theo1_cv_CDF_estimate_to_true_pointwise_ps} is equivalent to
\[
\forall \bx \in \cX, \forall s \in \left[ 0,1 \right],  \quad \sum_{j=1}^{n} w_j^b \left( \bx \right) \ind_{\left\lbrace U_j \leqslant s \right\rbrace} \overset{a.s.}{\underset{n \rightarrow \infty}{\longrightarrow}} s \ .
\]
As in the proof of Glivenko–Cantelli's Theorem, using that $s \mapsto \sum\limits_{j=1}^{n} w_j^b \left( \bx \right) \ind_{\left\lbrace U_j \left( \omega \right) \leqslant s \right\rbrace}$ is increasing and Dini's second theorem, we get the uniform convergence almost everywhere, which concludes the proof of the theorem. \hfill\qedsymbol

We now turn to the \veroc{proof of \Cref{prop:5:cv_diam_to_true} while that of \Cref{lem:5:cv_diam_to_bootstrap} is postponed to the Appendix.} To that aim, the following lemma, based on Vapnik-Chervonenkis classes \cite{vapnik71uniform} is a key tool. \veroc{The proof is postponed to the Appendix.}
\begin{lem}\label{lem:5:lemma_VC}
Consider $\cD_n$ and  $\cD_n^{\diamond}$, two independent datasets of independent $n$ samples of $(\bX,Y)$. Build a tree using $\cD_n$ with bootstrap and bagging procedure driven by $\Theta$. As before, $N^b \left( A_n \left( \Theta \right) \right) = N_n^b \left( \bx;\Theta, \cD_n \right)$ is the number of bootstrap observations of $\cD_n$ that fall into in $A_n \left( \Theta \right) = A_n \left( \bx;\Theta,\cD_n \right)$ and $N^{\diamond} \left( A_n \left( \Theta \right) \right) = N_n^{\diamond} \left( \bx;\Theta, \bX^{\diamond 1}, \ldots, \bX^{\diamond n}, \cD_n \right)$, the number of observations of $\cD_n^{\diamond}$ that fall into in $A_n \left( \Theta \right)$. Then,
\[
\forall \varepsilon > 0, \quad \PP \left( \left| N^b \left( A_n \left( \Theta \right) \right) - N^{\diamond} \left( A_n \left( \Theta \right) \right) \right| > \varepsilon \right) \leqslant 24 (n+1)^{2d}e^{- \varepsilon^2/288n} \ .
\]
\end{lem}

\Cref{lem:5:lemma_VC} is the main ingredient of the proof of \Cref{prop:5:cv_diam_to_true}.
\begin{myproof}[Proof of \Cref{prop:5:cv_diam_to_true}]
We aim to prove
\[
\forall \bx \in \cX, \forall y \in \RR, \quad \vero{\PP \left( F_{k,n}^{\diamond} \left( \left. y \right| \bX=\bx \right) \underset{n \rightarrow \infty}{\longrightarrow} F \left( \left. y \right| \bX=\bx \right) \right)} = 1 \ .
\]

Let $\bx \in \cX$ and $y \in \RR$, we have
\[
\left| F_{k,n}^{\diamond}  \left( \left. y \right| \bx \right) - F \left( \left. y \right| \bx \right) \right| \leqslant \left| \sum_{j=1}^{n} w_j^{\diamond} \left( \bx \right) \left( \ind_{\{ Y^{\diamond j} \leqslant y \}} -  F \left( \left. y \right| \bX^{\diamond j} \right) \right) \right| + \left| \sum_{j=1}^{n} w_j^{\diamond} \left( \bx \right) \left( F \left( \left. y \right| \bX^{\diamond j} \right) - F \left( \left. y \right| \bx \right) \right) \right| \ .
\]
\vero{Define 
\[
\displaystyle W_n =  \sum_{j=1}^{n} w_j^{\diamond} \left( \bx \right) \left( \ind_{\{ Y^{\diamond j} \leqslant y \}} -  F \left( \left. y \right| \bX^{\diamond j} \right) \right) = \sum_{j=1}^{n} w_j^{\diamond} \left( \bx \right) Z_j^{\diamond}
\]
with $Z_j^{\diamond} = \ind_{\{ Y^{\diamond j} \leqslant y \}} -  F \left( \left. y \right| \bX^{\diamond j} \right)$, $n$ i.i.d random variables and
\[
\displaystyle V_n = \sum_{j=1}^{n} w_j^{\diamond} \left( \bx \right) \left( F \left( \left. y \right| \bX^{\diamond j} \right) - F \left( \left. y \right| \bx \right) \right) \ .
\]
Remark that $\EE \left[ \left. Z_j^{\diamond} \right| \bX^{\diamond j} \right] = 0$.}\\
We first show that $\left( W_{n} \right)_{n \geqslant 1}$ goes to $0$ a.s. in the case of \Cref{hyp:4:num_2} item 2. This is achieved by adapting Hoeffding inequality's proof to our random weighted  sum context.
For any $\varepsilon >0$, $t \in \RR_{+}^*$, we have
\[
\PP \left( W_n > \varepsilon \right) \leq \EE \left[ e^{t W_n} \right] \cdot e^{-t\varepsilon} \ .
\]
We shall make use of the folklore lemma below.
\begin{lem*}
Let $X$ be a centred random variable, a.s. bounded by $1$. Then, for any \vero{$t \in \RR$}, $\EE \left[ e^{tX} \right] \leqslant e^{\frac{t^2}{2}}$.
\end{lem*}
Let $t>0$, we have
\begin{eqnarray*}
\EE \left[ e^{tW_n} \right] &=& \EE \left[ \prod_{j=1}^n e^{t  w_j^{\diamond}(\bx) Z_j^{\diamond}}\right] = \EE \left[ \EE \left[ \left. \prod_{j=1}^n e^{t  w_j^{\diamond}(\bx) Z_j^{\diamond}} \right| \cD_n, \Theta_1, \ldots, \Theta_k, \bX^{\diamond 1}, \ldots,\bX^{\diamond n} \right] \right] \\
\lefteqn{\mbox{conditionally to }  \cD_n, \Theta_1, \ldots, \Theta_k, \bX^{\diamond 1}, \ldots, \bX^{\diamond n}, \ \mbox{the $w_j^\diamond$ are constant and the $Z_j^\diamond$ are centred,}} \\
\lefteqn{\mbox{independent and bounded by $1$. Thus, using the folklore lemma,} } \\
\EE \left[ e^{tW_n} \right] & =&  \EE \left[ \prod_{j=1}^n \EE \left[ \left. e^{t  w_j^{\diamond}(\bx) Z_j^{\diamond}} \right| \cD_n, \Theta_1, \ldots, \Theta_k, \bX^{\diamond 1}, \ldots, \bX^{\diamond n} \right] \right] \leqslant \EE \left[ \prod_{j=1}^n e^{t^2 w_j^{\diamond}(\bx)^2 /2} \right].
\end{eqnarray*}
Let $K>0$ be such that for all $\ell=1\/, \ldots \/, k$, $N_n^b \left( A_n ( \ell ) \right) = N_n^b \left( \bx;\Theta_\ell,\cD_n \right) \geqslant K \sqrt{n} \left( \ln \left( n \right) \right)^\beta$ a.s. by using \Cref{hyp:4:num_2} item 2. Denote $\Gamma(\ell)$ the event $\left\lbrace N_n^\diamond(A_n(\ell)) < \frac{K \sqrt{n} \left( \ln (n) \right)^\beta}{2} \right\rbrace$. Remark that $\Gamma(\ell)\subset \left\lbrace \left| N_n^\diamond(A_n(\ell)) - N_n^b \left( A_n(\ell) \right) \right| >\frac{K \sqrt{n} \left( \ln (n) \right)^\beta}{2}\right\rbrace$. Thus, using \Cref{lem:5:lemma_VC}, we have that $\PP \left( \Gamma(\ell) \right) \leq 24 \left( n+1 \right)^{2d} \exp \left[ -\frac{K^2 \left( \ln (n) \right)^{2\beta}}{1152} \right]$.

We have
\begin{align*}
\sum_{j=1}^n  w_j^{\diamond}(\bx)^2 &= \sum_{j=1}^n \frac{w_j^{\diamond}(\bx)}{k} \left( \sum_{\ell=1}^k \frac{\ind_{\left\lbrace \bX^{\diamond j} \in A_n(\ell) \right\rbrace}}{N_n^\diamond(A_n(\ell))} 
\left( \ind_{\left\lbrace \Gamma(\ell)^c \right\rbrace} + \ind_{\left\lbrace \Gamma(\ell) \right\rbrace} \right) \right) \\
&\leqslant \sum_{j=1}^n w_j^{\diamond}(\bx)\left( \frac{2}{K \sqrt{n}(\ln n)^\beta} + \frac1k \sum_{\ell=1}^k \ind_{\left\lbrace \bX^{\diamond j} \in A_n(\ell) \right\rbrace} \ind_{\left\lbrace \Gamma(\ell) \right\rbrace} \right) \ .
\end{align*}

So that,
\begin{align*}
\EE \left[ \prod_{j=1}^n e^{t^2 w_j^{\diamond}(\bx)^2 /2} \right] &\leqslant \exp \left[ t^2 / \left( K \sqrt{n} \left( \ln (n) \right)^\beta \right) \right] \times \EE \left[ \exp{\left( \frac{t^2}{2} \cdot \ind_{\left\lbrace \bigcup_{\ell=1}^k \Gamma(\ell) \right\rbrace} \right)} \right] \\
&\leqslant \exp \left[ t^2 / \left( K \sqrt{n} \left( \ln (n) \right)^\beta \right) \right] \times \left( 1 + e^{t^2/2} \sum_{\ell=1}^k \PP \left( \Gamma(\ell) \right) \right) \\
&\leqslant \exp \left[ t^2 / \left( K \sqrt{n} \left( \ln (n) \right)^\beta \right) \right] \times \left( 1 + 24 k \left( n+1 \right)^{2d} \exp \left[\frac{t^2}{2} -\frac{K^2 \left( \ln (n) \right)^{2\beta}}{1152} \right] \right) \ .
\end{align*}
Taking $t^2 = \frac{K^2 \left( \ln (n) \right)^{2\beta}}{576}$ leads to
\[
\PP \left( W_n > \varepsilon \right) \leqslant \left( 1 + 24 k \left( n+1 \right)^{2d} \right) \exp \left[\frac{K \left( \ln (n) \right)^{\beta}}{576 \sqrt{n}} -\frac{\varepsilon K \left( \ln (n) \right)^{\beta}}{24} \right] \ .
\]
The same upper bound is obtained for $\PP \left( W_n < -\varepsilon \right)$ by using that $\PP \left( W_n < -\varepsilon \right) = \PP \left( -W_n > \varepsilon \right)$.

Thus, by using \Cref{hyp:4:num_2}, item 1., $k=O(n^\alpha)$ so that the right hand side is summable,  we conclude that $W_n$ goes to $0$ almost surely.

\vero{In the case where  \Cref{hyp:4:num_2} item 3 is satisfied, the proof that $\left( W_{n} \right)_{n \geqslant 1}$  goes to $0$ a.s. is done in a similar spirit, \veroc{the sketch of proof under \Cref{hyp:4:num_2} item 3, is proposed in the Appendix. }

Finally, we show that $\left( V_{n} \right)_{n \geqslant 1}$ goes to 0 a.s. which easily follows from Assumption \ref{hyp:4:num_1}.}
This allows us to conclude that
\[
\forall \bx \in \cX, \forall y \in \RR, \quad F_{k,n}^{\diamond} \left( \left. y \right| \bX=\bx \right) \overset{a.s.}{\underset{n \rightarrow \infty}{\longrightarrow}} F \left( \left. y \right| \bX=\bx \right) \ .
\]
\end{myproof}

\veroc{Thus, we have proved  \Cref{theo:4:num_1}, the proof of \Cref{theo:4:num_2} which is a bit simpler is left to the reader.}

\section{Conclusion}\label{sec:7:conclusion}

This article proposes two conditional distribution functions and conditional quantiles approximations based on random forests. The former is a natural generalisation of the random forest estimator of the regression function making use of the bootstrap samples, while the latter is based on a variant using only the original dataset. 

The consistency of the bootstrap samples based estimator is shown under realistic assumptions and constitutes the major contribution of this paper. Indeed, this is the first consistency result handling the bootstrap component in a random forest method whereas it is usually replaced by subsampling. As for the second estimator, the consistency proof established in \cite{meinshausen2006quantile} for a simplified random forest model is extended to a realistic one by taking into account all the randomness used in the procedure. The two estimators have close performances on our toy example. A specific interest of the bootsrap estimation is that the Out-Of-Bag samples could be used for cross-validation and / or back-testing procedures.

The estimators developed in this paper rest on trees grown with the CART-split criterion. But the assumptions providing the consistency results are detached from the split procedure used. Thus, the theoretical tools developed here could be useful for a large class of methods by just changing the splitting scheme. An ambitious additional work would be to develop a theoretical analysis for obtaining convergence rates and also to construct confidence intervals.

 \section*{Acknowledgments}
 We are grateful to Andr\'es Cuberos, Ecaterina Nisipasu, Mathieu Poulin and Przemyslaw Sloma from SCOR for their valuable comments and support. We are also much indebted to Roland Denis and Benoit Fabrèges for intensive support on computational aspects. We are grateful to anonymous reviewers of  previous versions of the paper, their comments helped to improve the paper.



\bibliographystyle{plain}

\bibliography{bibliography}

\begin{thebibliography}{10}

\bibitem{amit1997shape}
Yali Amit and Donald Geman.
\newblock Shape quantization and recognition with randomized trees.
\newblock {\em Neural computation}, 9(7):1545--1588, 1997.

\bibitem{bezanson2017julia}
Jeff Bezanson, Alan Edelman, Stefan Karpinski, and Viral~B Shah.
\newblock Julia: A fresh approach to numerical computing.
\newblock {\em SIAM review}, 59(1):65--98, 2017.

\bibitem{biau2010layered}
G{\'e}rard Biau and Luc Devroye.
\newblock On the layered nearest neighbour estimate, the bagged nearest
  neighbour estimate and the random forest method in regression and
  classification.
\newblock {\em Journal of Multivariate Analysis}, 101(10):2499--2518, 2010.

\bibitem{biau2012analysis}
Gérard Biau.
\newblock Analysis of a random forests model.
\newblock {\em Journal of Machine Learning Research}, 13(Apr):1063--1095, 2012.

\bibitem{breiman1996bagging}
Leo Breiman.
\newblock Bagging predictors.
\newblock {\em Machine learning}, 24(2):123--140, 1996.

\bibitem{breiman2001random}
Leo Breiman.
\newblock Random forests.
\newblock {\em Machine learning}, 45(1):5--32, 2001.

\bibitem{breiman2004consistency}
Leo Breiman.
\newblock Consistency for a simple model of random forests.
\newblock 2004.

\bibitem{breiman1984classification}
Leo Breiman, Jerome~H Friedman, Richard~A Olshen, and Charles~J Stone.
\newblock Classification and regression trees.
\newblock 1984.

\bibitem{browne2017estimate}
Thomas Browne, Jean-Claude Fort, Bertrand Iooss, and Lo{\"\i}c Le~Gratiet.
\newblock Estimate of quantile-oriented sensitivity indices.
\newblock 2017.

\bibitem{devroye2013probabilistic}
Luc Devroye, L{\'a}szl{\'o} Gy{\"o}rfi, and G{\'a}bor Lugosi.
\newblock {\em A probabilistic theory of pattern recognition}, volume~31.
\newblock Springer Science \& Business Media, 2013.

\bibitem{dietterich2000ensemble}
Thomas~G Dietterich.
\newblock Ensemble methods in machine learning.
\newblock In {\em International workshop on multiple classifier systems}, pages
  1--15. Springer, 2000.

\bibitem{Efron_1979}
B.~Efron.
\newblock Bootstrap methods: Another look at the jackknife.
\newblock {\em The Annals of Statistics}, 7:1--26, 1979.

\bibitem{qosa_indices}
K{\'e}vin Elie-Dit-Cosaque.
\newblock qosa-indices, a python package available at: {\tt
  https://gitlab.com/qosa\_index/qosa}, 2020.

\bibitem{qosa_shapley}
K{\'e}vin Elie-Dit-Cosaque and V{\'e}ronique Maume-Deschamps.
\newblock Goal-oriented shapley effects with a special attention to the
  quantile-oriented case.
\newblock {\em SIAM/ASA Journal on Uncertainty Quantification - JUQ}, to
  appear.

\bibitem{ConditionalDistributionForest}
Benoit Fabr{\`e}ge and V{\'e}ronique Maume-Deschamps.
\newblock Conditional distribution forest: a julia package available at {\tt
  https://github.com/bfabreges/conditionaldistributionforest.jl}, 2020.

\bibitem{fort2016new}
Jean-Claude Fort, Thierry Klein, and Nabil Rachdi.
\newblock New sensitivity analysis subordinated to a contrast.
\newblock {\em Communications in Statistics-Theory and Methods},
  45(15):4349--4364, 2016.

\bibitem{goehry2019random}
Benjamin Goehry.
\newblock Random forests for time-dependent processes.
\newblock 2019.

\bibitem{gyorfi2006distribution}
L{\'a}szl{\'o} Gy{\"o}rfi, Michael Kohler, Adam Krzyzak, and Harro Walk.
\newblock {\em A distribution-free theory of nonparametric regression}.
\newblock Springer Science \& Business Media, 2006.

\bibitem{ho1998random}
Tin~Kam Ho.
\newblock The random subspace method for constructing decision forests.
\newblock {\em IEEE transactions on pattern analysis and machine intelligence},
  20(8):832--844, 1998.

\bibitem{klusowki}
Jason~M. Klusowski.
\newblock Analyzing cart.
\newblock 2020.

\bibitem{koenker2001quantile}
Roger Koenker and Kevin~F Hallock.
\newblock Quantile regression.
\newblock {\em Journal of economic perspectives}, 15(4):143--156, 2001.

\bibitem{lin2006random}
Yi~Lin and Yongho Jeon.
\newblock Random forests and adaptive nearest neighbors.
\newblock {\em Journal of the American Statistical Association},
  101(474):578--590, 2006.

\bibitem{maume2018estimation}
V{\'e}ronique Maume-Deschamps and Ibrahima Niang.
\newblock Estimation of quantile oriented sensitivity indices.
\newblock {\em Statistics \& Probability Letters}, 134:122--127, 2018.

\bibitem{maume2017quantile}
V{\'e}ronique Maume-Deschamps, Didier Rulli{\`e}re, and A~Usseglio-Carleve.
\newblock Quantile predictions for elliptical random fields.
\newblock {\em Journal of Multivariate Analysis}, 159:1--17, 2017.

\bibitem{meinshausen2006quantile}
Nicolai Meinshausen.
\newblock Quantile regression forests.
\newblock {\em Journal of Machine Learning Research}, 7(Jun):983--999, 2006.

\bibitem{qrf}
Nicolai Meinshausen.
\newblock Quantile regression forests, a r package available at {\tt
  https://cran.r-project.org/package=quantregforest}., 2019.

\bibitem{mentch2016quantifying}
Lucas Mentch and Giles Hooker.
\newblock Quantifying uncertainty in random forests via confidence intervals
  and hypothesis tests.
\newblock {\em The Journal of Machine Learning Research}, 17(1):841--881, 2016.

\bibitem{radial_functions}
Jooyoung Park and Irwin~W Sandberg.
\newblock Universal approximation using radial-basis-function networks.
\newblock {\em Neural computation}, 1991.

\bibitem{scikit-learn}
F.~Pedregosa, G.~Varoquaux, A.~Gramfort, V.~Michel, B.~Thirion, O.~Grisel,
  M.~Blondel, P.~Prettenhofer, R.~Weiss, V.~Dubourg, J.~Vanderplas, A.~Passos,
  D.~Cournapeau, M.~Brucher, M.~Perrot, and E.~Duchesnay.
\newblock Scikit-learn: Machine learning in {P}ython.
\newblock {\em Journal of Machine Learning Research}, 12:2825--2830, 2011.

\bibitem{Rsoftware}
{R Core Team}.
\newblock {\em R: A Language and Environment for Statistical Computing}.
\newblock R Foundation for Statistical Computing, Vienna, Austria, 2019.

\bibitem{scornet2016asymptotics}
Erwan Scornet.
\newblock On the asymptotics of random forests.
\newblock {\em Journal of Multivariate Analysis}, 146:72--83, 2016.

\bibitem{scornetpromenade}
Erwan Scornet.
\newblock Promenade en for{\^e}ts al{\'e}atoires.
\newblock {\em MATAPLI}, 111, 2016.

\bibitem{scornet2016random}
Erwan Scornet.
\newblock Random forests and kernel methods.
\newblock {\em IEEE Transactions on Information Theory}, 62(3):1485--1500,
  2016.

\bibitem{scornet2015supplementary}
Erwan Scornet, G{\'e}rard Biau, and Jean-Philippe Vert.
\newblock Supplementary materials for: Consistency of random forests.
\newblock {\em arXiv}, 1510, 2015.

\bibitem{scornet2015consistency}
Erwan Scornet, G{\'e}rard Biau, Jean-Philippe Vert, et~al.
\newblock Consistency of random forests.
\newblock {\em The Annals of Statistics}, 43(4):1716--1741, 2015.

\bibitem{vapnik71uniform}
V.~N. Vapnik and A.~Ya. Chervonenkis.
\newblock On the uniform convergence of relative frequencies of events to their
  probabilities.
\newblock {\em Theory of Probability and its Applications}, 16(2):264--280,
  1971.

\bibitem{wager2018estimation}
Stefan Wager and Susan Athey.
\newblock Estimation and inference of heterogeneous treatment effects using
  random forests.
\newblock {\em Journal of the American Statistical Association},
  113(523):1228--1242, 2018.

\bibitem{wager2015adaptive}
Stefan Wager and Guenther Walther.
\newblock Adaptive concentration of regression trees, with application to
  random forests.
\newblock {\em arXiv preprint arXiv:1503.06388}, 2015.

\end{thebibliography}

\begin{appendix}
\section{Technical complements}\label{sec:proof2}

\vero{We provide some technical Lemmas used in Section \ref{sec:variation} and Section \ref{sec:5:proofs}. 
\begin{lem}\label{lem:vc_app}
Let $\cA$ be the class of rectangles in $\RR^d$, let $\kappa>0$, for any $n$,
\[
\PP \left( \sup_{B\in\cA} \left| \dfrac{N_B^b}n - \PP(\bX\in B)\right| > \kappa\right) \leq 16 \kev{(n+1)}^{2d}e^{-n\kappa^2/\kev{128}} \ .
\]
\end{lem}
\begin{proof}
Denote by $N_B$ the number of elements of the original sample $\cD_n$ that are in $B$. We have:
\begin{eqnarray*}
 \lefteqn{\PP\left(\sup_{B\in\cA}\left|\frac{N_B^b}n - \PP(\bX\in B)\right|>\kappa\right) }\\
 &\leq & \PP\left(\sup_{B\in\cA}\left|\frac{N_B}n - \PP(\bX\in B)\right|>\kappa/2 \right)+ \PP\left(\sup_{B\in\cA}\left|\frac{N_B^b}n - \frac{N_B}n\right|>\kappa/2\right)\\
 &=&\PP\left(\sup_{B\in\cA}\left|\frac{N_B}n - \PP(\bX\in B)\right|>\kappa/2 \right) + \EE\left(\PP \left( \left. \sup_{B\in\cA}\left| \frac1n \sum_{i=1}^n \ind_{\left\lbrace \bX^{\star i} \in B \right\rbrace} - \PP(\bX^\star \in B | \cD_n)\right| >\kappa/2 \right| \cD_n\right)\right)
\end{eqnarray*}
where $\bX^\star$ is \veroc{a bootstrap sample}, whose distribution conditionally to $\cD_n$ is uniform on $\{\bX^1,\ldots,\bX^n\}$. We apply Vapnik-Chervonenkis (\cite{vapnik71uniform}) inegality conditionally to $\cD_n$ to get that the second term is bounded above by $8 \kev{(n+1)}^{2d}e^{-n\kappa^2/\kev{128}}$ and  Vapnik-Chervonenkis inegality to get that the first term admits the same above bound. 
\end{proof}
\begin{lem}\label{lem:giorfy}
 Let $(\bX_i^\star\/,Y_i^\star)$ be \veroc{a bootstrap sample} from $(\bX_1\/,Y_1)\/, \ldots \/, (\bX_n\/,Y_n)$.  For any \verob{$D>0$}, for any $\kappa>0$, $A\in \cA$ and $\frac1p+\frac1q=1$. Assume that $\EE(|Y|^p)<\infty$ then
\verob{\begin{eqnarray*}
\lefteqn{\PP\left(\left|\frac1n \sum_{i=1}^nY_i^\star \ind_{\{\bX_i^\star\in A\}} - \EE(Y\ind_{\{\bX \in A\}})\right| >\kappa \right)}  \\
 &\leq& 6\left(\frac{24eD}\kappa \ln \left(\frac{48eD}\kappa\right)\right)^{2d} e^{-n\kappa^2/(128 D^2)} + C\frac{\EE(Y^p)^\frac1p \PP(Y>D)^\frac1q}{\kappa}\/.
 \end{eqnarray*}
 }
\end{lem}
\begin{proof}
\verob{Let $Z=\min(Y\/,D)$, $Z_i=\min(Y_i\/,D)$, $Z_i^\star=\min(Y_i^\star\/,D)$.} Write:
\begin{eqnarray*}
\lefteqn{\left|\frac1n \sum_{i=1}^nY_i^\star \ind_{\{\bX_i^\star\in A\}} - \EE(Y\ind_{\{\bX \in A\}})\right| }\\
&\leq & \EE((Y-Z)\ind_{\{\bX \in A\}}) + \left|\frac1n\sum_{i=1}^n (Y_i^\star -Z_i^\star)\ind_{\{\bX_i^\star \in A\}}\right|+ \\
&&  \left|\frac1n \sum_{i=1}^n Z_i^\star \ind_{\{\bX_i^\star \in A\}} - \EE(Z\ind_{\{\bX\in A\}})\right| \/.
\end{eqnarray*}
The probability that the sum of the first two terms is greater than $\kappa$ is bounded above by \verob{$C\frac{\EE(Y^p)^\frac1p \PP(Y>D)^\frac1q}{\kappa}$}. The probability that the \kev{last} term is greater than $\kappa$ is bounded by following the proof of Theorem 9.6 in \cite{gyorfi2006distribution} page 155, once conditionally to $\cD_n$ to take into account the bootstrap sample and then unconditionally.
\end{proof}

}

\veroc{

Let us now turn to the proof of \Cref{lem:5:cv_diam_to_bootstrap} which shows that the dummy estimator $F_{k,n}^{\diamond}$ is close to the interesting one $F_{k,n}^b$.
\begin{myproof}[Proof of \Cref{lem:5:cv_diam_to_bootstrap}]
Let $\bx \in \cX$ and $y \in \RR$, recall that $B_j\left( \Theta_\ell^1,\cD_n \right)$ has been defined in Section \ref{sub_sec:3.1:QRF}.We have
\begin{align*}
& \left| F_{k,n}^{\diamond} \left( \left. y \right| \bX=\bx \right) - F_{k,n}^b \left( \left. y \right| \bX=\bx \right) \right| \\
&= \left| \dfrac{1}{k} \sum_{\ell=1}^k \left( \sum_{j=1}^{n} \dfrac{\ind_{\left\lbrace \bX^{\diamond j} \in A_n \left( \bx;\Theta_\ell,\cD_n \right) \right\rbrace} \ind_{\left\lbrace Y^{\diamond j} \leqslant y \right\rbrace}}{N_n^{\diamond} \left( \bx; \Theta_\ell,\bX^{\diamond 1}, \ldots, \bX^{\diamond n}, \cD_n \right)} - \sum_{j=1}^{n} \dfrac{B_{j} \left( \Theta_\ell^1,\cD_n \right) \ind_{\left\lbrace \bX^j \in A_n \left( \bx;\Theta_\ell,\cD_n \right) \right\rbrace} \ind_{\left\lbrace Y^{j} \leqslant y \right\rbrace}}{N_n^b \left( \bx; \Theta_\ell, \cD_n \right)} \right) \right| \\
&= \left| \dfrac{1}{k} \sum_{\ell=1}^k \left( \dfrac{\# \left\lbrace j \leqslant J^{\diamond} \ / \ \bX^{\diamond (j)} \in A_n \left( \Theta_\ell \right) \right\rbrace}{N^{\diamond} \left( A_n \left( \Theta_\ell \right) \right)} - \dfrac{\sum\limits_{j \in \cS} B_{j} \left( \Theta_\ell^1,\cD_n \right)}{N^b \left( A_n \left( \Theta_\ell \right) \right)} \right) \right| \textnormal{ with } \cS = \left\lbrace j \leqslant J \ / \ \bX^{(j)} \in A_n \left( \Theta_\ell \right) \right\rbrace \ ,
\end{align*}
where we denote $A_n \left( \Theta_\ell \right) = A_n \left( \bx;\Theta_\ell,\cD_n \right), N^{\diamond} \left( A_n \left( \Theta_\ell \right) \right) = N_n^{\diamond} \left( \bx; \Theta_\ell,\bX^{\diamond 1}, \ldots, \bX^{\diamond n}, \cD_n \right)$ and $N^b \left( A_n \left( \Theta_\ell \right) \right) = N_n^b \left( \bx; \Theta_\ell, \cD_n \right)$. $J, J^\diamond$ are such that $Y^{\diamond \left( J^\diamond \right)} \leqslant y < Y^{\diamond \left( J^{\diamond}+1 \right)}$ and $Y^{\left( J \right)} \leqslant y < Y^{\left( J+1 \right)}$, with $Y^{\diamond (j)}$ (resp. $Y^{(j)}$) the order statistics of $(Y^{\diamond 1}, \ldots, Y^{\diamond n})$ (resp. $(Y^1, \ldots, Y^n)$) and the $\bX^{\diamond (j)}$  (resp. $\bX^{(j)}$) the corresponding $\bX^{\diamond p}$'s (resp. $\bX^p$'s).

Let us consider for some $\ell \in \llbracket 1,k \rrbracket$,
\[
G = \dfrac{\# \left\lbrace j \leqslant J^{\diamond} \ / \ \bX^{\diamond (j)} \in A_n \left( \Theta_\ell \right) \right\rbrace}{N^{\diamond} \left( A_n \left( \Theta_\ell \right) \right)} - \dfrac{\sum\limits_{j \in \cS} B_{j} \left( \Theta_\ell^1,\cD_n \right)}{N^b \left( A_n \left( \Theta_\ell \right) \right)} \overset{\textnormal{def}}{=} \dfrac{N_{J^{\diamond}}^{\diamond} \left( A_n \left( \Theta_\ell \right) \right)}{N^{\diamond} \left( A_n \left( \Theta_\ell \right) \right)} - \dfrac{N_J \left( A_n \left( \Theta_\ell \right) \right)}{N^b \left( A_n \left( \Theta_\ell \right) \right)} \ .
\]
We have, 
\begin{align*}
\left| G \right| &\leqslant \dfrac{\left| N^{\diamond} \left( A_n \left( \Theta_\ell \right) \right) - N^b \left( A_n \left( \Theta_\ell \right) \right) \right|}{N^b \left( A_n \left( \Theta_\ell \right) \right)} + \dfrac{\left| N_{J^{\diamond}}^{\diamond} \left( A_n \left( \Theta_\ell \right) \right) - N_J \left( A_n \left( \Theta_\ell \right) \right) \right|}{N^b \left( A_n \left( \Theta_\ell \right) \right)} \\
&\overset{\textnormal{def}}{=} \left| G_1 \right| + \left| G_2 \right|
\end{align*}

We continue the proof below in the case where \Cref{hyp:4:num_2} item 3. is satisfied. The case where item 2. is verified is done easier following the same lines. Let $\varepsilon >0$. \vero{We are now going to show the almost everywhere convergence to 0 for each term $G_1$ and $G_2$. Let us start with $G_1$.}
\begin{align*}
\PP \left( \left| G_1 \right| > \varepsilon \right) &= \PP \left( \dfrac{\left| N^{\diamond} \left( A_n \left( \Theta_\ell \right) \right) - N^b \left( A_n \left( \Theta_\ell \right) \right) \right|}{N^b \left( A_n \left( \Theta_\ell \right) \right)} > \varepsilon \right) \\
&= \PP \left( \left| N^{\diamond} \left( A_n \left( \Theta_\ell \right) \right) - N^b \left( A_n \left( \Theta_\ell \right) \right) \right| > \varepsilon N^b \left( A_n \left( \Theta_\ell \right) \right), N^b \left( A_n \left( \Theta_\ell \right) \right) > \lambda \right) \\
&\quad + \PP \left( \left| N^{\diamond} \left( A_n \left( \Theta_\ell \right) \right) - N^b \left( A_n \left( \Theta_\ell \right) \right) \right| > \varepsilon N^b \left( A_n \left( \Theta_\ell \right) \right), N^b \left( A_n \left( \Theta_\ell \right) \right) \leqslant \lambda \right) \\ \intertext{ where $\lambda = \dfrac{\EE \left[ N^b \left( A_n \left( \Theta \right) \right) \right]}{2}$} 
&\leqslant \PP \left( \left| N^{\diamond} \left( A_n \left( \Theta_\ell \right) \right) - N^b \left( A_n \left( \Theta_\ell \right) \right) \right| > \varepsilon \lambda \right) + \PP \left( N^b \left( A_n \left( \Theta_\ell \right) \right) \leqslant \lambda \right)
\end{align*}

\verob{Thanks to} Bienaymé-Tchebychev's inequality,
\begin{align}
\PP \left( N^b \left( A_n \left( \Theta_\ell \right) \right) \leqslant \lambda \right) &\leqslant 4 \dfrac{\var \left( N^b \left( A_n \left( \Theta_\ell \right) \right) \right)}{\left( \EE \left[ N^b \left( A_n \left( \Theta_\ell \right) \right) \right] \right)^2} \nonumber \\
&\leqslant 4 \left( \CV \left( N^b \left( A_n \left( \Theta \right) \right) \right) \right)^2 \ . \label{eq:5:upper_bound_with_CV}
\end{align}
Now, using \Cref{lem:5:lemma_VC} and \Cref{hyp:4:num_2}, we get
\begin{align*}
\PP \left( \left| G_1 \right| > \varepsilon \right) &\leqslant \PP \left( \left| N^{\diamond} \left( A_n \left( \Theta_\ell \right) \right) - N^b \left( A_n \left( \Theta_\ell \right) \right) \right| > \varepsilon \lambda \right) + 4 \left( \CV \left( N^b \left( A_n \left( \Theta \right) \right) \right) \right)^2 \\
&\leqslant 24(n+1)^{2d} \exp \left[- \dfrac{\varepsilon^2 K^2 \left( \ln \left( n \right) \right)^{2 \beta}}{1152} \right] + \dfrac{4M^2}{n^{\kev{\alpha + 1}}  \left( \ln \left( n \right) \right)^\gamma} \ .
\end{align*}
Then, thanks to Borel–Cantelli Lemma
\[
\forall \varepsilon > 0, \quad \PP \left( \limsup\limits_{n \rightarrow \infty} \left\lbrace \left| G_1 \right| > \varepsilon \right\rbrace \right) = 0 \ ,
\]
which implies $\displaystyle G_1 \overset{a.s.}{\underset{n \rightarrow \infty}{\longrightarrow}} 0$.\\
\vero{The $G_2$ term is treated by using 
%
%
again the Vapnik-Chervonenkis theory. By considering the class $\cB = \left\lbrace \prod\limits_{i=1}^{d} \left[ a_i, b_i \right] \times \left] -\infty, y \right] : a_i, b_i \in \overline{\RR} \right\rbrace$, it gives (following the lines of proof of Lemma \ref{lem:vc_app})
\begin{align*}
\PP \left( \left| G_2 \right| > \varepsilon \right)  &\leqslant \PP \left( \left| N_{J^{\diamond}}^{\diamond} \left( A_n \left( \Theta_\ell \right) \right) - N_J \left( A_n \left( \Theta_\ell \right) \right) \right| > \varepsilon \lambda \right) + 4 \left( \CV \left( N \left( A_n \left( \Theta \right) \right) \right) \right)^2 \\
&\leqslant 24(n+1)^{2d} \exp \left[- \dfrac{\varepsilon^2 K^2 \left( \ln \left( n \right) \right)^{2 \beta}}{1152} \right] + \dfrac{4M^2}{n^{\kev{\alpha + 1}}  \left( \ln \left( n \right) \right)^\gamma} \ .
\end{align*}
}
Thanks to Borel–Cantelli Lemma, we get
\[
\forall \varepsilon > 0, \quad \PP \left( \limsup\limits_{n \rightarrow \infty} \left\lbrace \left| G_2 \right| > \varepsilon \right\rbrace \right) = 0 \ ,
\]
which implies that $\displaystyle G_2 \overset{a.s.}{\underset{n \rightarrow \infty}{\longrightarrow}} 0$.

We conclude that $G$ goes to $0$ for all $\ell$, thus,
\[
\forall \bx \in \cX, \forall y \in \RR, \quad \left| F_{k,n}^{\diamond} \left( \left. y \right| \bX=\bx \right) -  F_{k,n}^b \left( \left. y \right| \bX=\bx \right) \right| \overset{a.s.}{\underset{n \rightarrow \infty}{\longrightarrow}} 0 \ .
\]

In the case where \Cref{hyp:4:num_2} item 2. is verified, it exists $K>0$ such that
\[
N^b \left( A_n \left( \Theta_\ell \right) \right) \geqslant K \sqrt{n} \left( \ln \left( n \right) \right)^\beta \ a.s.
\]
So that $\PP(|G_1|>\varepsilon)$ and $\PP(|G_2|>\varepsilon)$ are bounded above respectively by 
\begin{itemize}
\item $ \PP \left( \left| N^b \left( A_n \left( \Theta_\ell \right) \right) - N^{\diamond} \left( A_n \left( \Theta_\ell \right) \right) \right| > \varepsilon K \sqrt{n} \left( \ln \left( n \right) \right)^\beta \right)$,
\item and $\PP \left( \left|  N_{J^{\diamond}}^{\diamond} \left( A_n \left( \Theta_\ell \right) \right) - N_J \left( A_n \left( \Theta_\ell \right) \right) \right| > \varepsilon K \sqrt{n} \left( \ln \left( n \right) \right)^\beta \right)$.
\end{itemize}
A simple application of \Cref{lem:5:lemma_VC} and an adaptation of it to $N_J \left( A_n \left( \Theta_\ell \right) \right)$ show that $G_1$ and $G_2$ go to $0$ a.s.
\end{myproof}

Now we turn to \Cref{lem:5:lemma_VC} which is key for the proof of Proposition \ref{prop:5:cv_diam_to_true} and that of \Cref{lem:5:cv_diam_to_bootstrap}. 
\begin{myproof}[Proof of \Cref{lem:5:lemma_VC}]
Let $\varepsilon > 0$ and $\bx \in \cX$, we have
\begin{align*}
& \PP \left( \left| N^b \left( A_n \left( \Theta \right) \right) - N^{\diamond} \left( A_n \left( \Theta \right) \right) \right| > \varepsilon \right) \\ 
&\leqslant \PP \left( \left| \dfrac{N^b \left( A_n \left( \Theta \right) \right)}{n} - \dfrac{1}{n} \sum_{j=1}^n \ind_{\left\lbrace \bX^j \in A_n \left( \Theta \right) \right\rbrace} \right| > \dfrac{\varepsilon}{3n} \right) + \PP \left( \left| \dfrac{1}{n} \sum_{j=1}^n \ind_{\left\lbrace \bX^j \in A_n \left( \Theta \right) \right\rbrace} - \PP_{\bX} \left( \bX \in A_n \left( \Theta \right) \right) \right| > \dfrac{\varepsilon}{3n} \right) \\
& \quad + \PP \left( \left| \dfrac{N^{\diamond} \left( A_n \left( \Theta \right) \right)}{n} - \PP_{\bX} \left( \bX \in A_n \left( \Theta \right) \right) \right| > \dfrac{\varepsilon}{3n} \right) \\
&\leqslant \PP \left( \sup_{A \in \cB} \left| \dfrac{1}{n} \sum_{j=1}^n  \ind_{\left\lbrace \bX^{*j} \in A \right\rbrace} - \dfrac{1}{n} \sum_{j=1}^n \ind_{\left\lbrace \bX^j \in A \right\rbrace} \right| > \dfrac{\varepsilon}{3n} \right) \\
& \quad + \PP \left( \sup_{A \in \cB} \left| \dfrac{1}{n} \sum_{j=1}^n \ind_{\left\lbrace \bX^j \in A \right\rbrace} - \PP_{\bX} \left( \bX \in A \right) \right| >  \dfrac{\varepsilon}{3n} \right) + \PP \left( \sup_{A \in \cB} \left| \dfrac{1}{n} \sum_{j=1}^n \ind_{\left\lbrace \bX^{\diamond j} \in A \right\rbrace} - \PP_{\bX} \left( \bX \in A \right) \right| >  \dfrac{\varepsilon}{3n} \right)
\end{align*}
where $\bX^{*1}\/, \ldots\/, \bX^{*n}$ denotes a bootstrap sample and $\cB = \left\lbrace \prod\limits_{i=1}^d \left[ a_i, b_i \right] : a_i, b_i \in \overline{\RR} \right\rbrace$. The last two right-hand terms are handled thanks to a direct application of the Theorem of \cite{vapnik71uniform} over the class $\cB$ whose Vapnik-Chervonenkis dimension is $2d$. This class is nothing more than an extension of the class $\cR$ of rectangles in $\RR^d$. Following the lines of the proof of Theorem 13.8 in \cite{devroye2013probabilistic}, one sees that the classes $\cR$ and $\cB$ have the same Vapnik-Chervonenkis dimension.

\vero{The first right hand term is treated by applying Vapnik-Chervonenkis' Theorem under the conditional distribution given $\cD_n$ as in the proof of Lemma \ref{lem:vc_app}.} 
%
%
%

Finally, we get the overall upper bound
\[
\PP \left( \left| N^b \left( A_n \left( \Theta \right) \right) - N^{\diamond} \left( A_n \left( \Theta \right) \right) \right| > \varepsilon \right) \leqslant 24 (n+1)^{2d}e^{- \varepsilon^2/288n} \ .
\]
\end{myproof}
Let us now sketch the part of proof  of Proposition \ref{prop:5:cv_diam_to_true} which states  that $\left( W_{n} \right)_{n \geqslant 1}$  goes to $0$ a.s., under \Cref{hyp:4:num_2} item 3.
\begin{proof}
\begin{enumerate}
 \item  First show that $\left( W_{n^2} \right)_{n \geqslant 1}$ goes to 0 a.s. This is achieved by decomposing 
 \begin{align*}
\EE \left[ \left( W_n \right)^2 \right] &= \EE \left[ \left( \sum_{j=1}^{n} w_j^{\diamond} \left( \bx \right) Z_j^{\diamond} \right)^2 \right] \\
&= \sum_{j=1}^{n} \sum_{m=1}^{n} \EE \left[ w_j^{\diamond} \left( \bx \right) w_m^{\diamond} \left( \bx \right) Z_j^{\diamond} Z_m^{\diamond} \right] \\
&= \sum_{j=1}^{n} \EE \left[ w_j^{\diamond^2} \left( \bx \right) Z_j^{\diamond^2} \right] + \sum\limits_{\substack{1 \leqslant j,m \leqslant n \\ j \neq m}} \EE \left[ w_j^{\diamond} \left( \bx \right) w_m^{\diamond} \left( \bx \right) Z_j^{\diamond} Z_m^{\diamond} \right] \\
&\overset{\textnormal{def}}{=} I_n + J_n  
\end{align*}
Bienaimé-Tchebychev's inequality, \Cref{lem:5:lemma_VC} and \Cref{hyp:4:num_2} items 1. and 3., give that there exist $C, K$ and $M$ positive constants such that
\begin{align*}
I_n &\leqslant k \PP \left( \left| N^b \left( A_n \left( \Theta \right) \right) - N^{\diamond} \left( A_n \left( \Theta \right) \right) \right| > \lambda \right) + \dfrac{4}{\EE \left[ N^b \left( A_n \left( \Theta \right) \right) \right]} + 4 \kev{k} \left( \CV \left( N^b \left( A_n \left( \Theta \right) \right) \right) \right)^2 \\
&\leqslant 24 C n^\alpha (n+1)^{2d} \exp \left[- \dfrac{K^2 \left( \ln \left( n \right) \right)^{2 \beta}}{4608} \right] +\dfrac{4}{K \sqrt{n} \left( \ln \left( n \right) \right)^\beta} + \dfrac{4 CM^2}{n  \left( \ln \left( n \right) \right)^\gamma} \ .
\end{align*}
Then, the trick of using a second sample $\cD_n^{\diamond}$ independent of the first-one and the random variable $\Theta$ is really important to handle the $J_n$ term. Indeed, we have $J_n = 0$ because $\EE \left[ \left. Z_m^{\diamond} \right| \bX^{\diamond m} \right] = 0$ while the equivalent term encountered in the proof of the Theorem 2 developed by \cite{scornet2015consistency} is handled using a conjecture regarding the correlation behavior of the CART algorithm that is difficult to verify (cf. assumption (H2) of \cite{scornet2015consistency}). 
Finally,
\[
\forall \varepsilon > 0, \quad \PP \left( \left| W_n \right| \geqslant \varepsilon \right) \leqslant \dfrac{\EE \left[ \left( W_n \right)^2 \right]}{\varepsilon^2} = \dfrac{I_n}{\varepsilon^2} \ .
\]
Hence, since $\displaystyle\sum_{n\geq 1} I_{n^2}<\infty$, Borel–Cantelli Lemma gives
\[
\forall \varepsilon > 0, \quad \PP \left( \limsup\limits_{n \rightarrow \infty} \left\lbrace \left| W_{n^2} \right| \geqslant \varepsilon \right\rbrace \right) = 0 \ ,
\]
which implies that $\displaystyle W_{n^2} \overset{a.s.}{\underset{n \rightarrow \infty}{\longrightarrow}} 0$.

\item Show that $\left( W_{n} \right)_{n \geqslant 1}$ converges almost surely to 0. Bienaimé-Tchebytchev inequality, \Cref{lem:5:lemma_VC} and \Cref{hyp:4:num_2} items 1. and 3. as well as Borel-Cantelli Lemma give that 
\[
\forall \varepsilon > 0, \quad \PP \left( \limsup\limits_{n \rightarrow \infty} \left\lbrace \left| W_{n} - W_{p^2} \right| \geqslant \varepsilon \right\rbrace \right) = 0 \ ,
\]
for $p=p \left( n \right) = \left\lfloor \sqrt{n} \right\rfloor$.
From this, we deduce that $\left( W_n \right)_{n \geqslant 1}$ goes to 0 a.s.
\end{enumerate}
\end{proof}

We end this Appendix with the proof that sums of products of functions considered in Remark \ref{rem:pique} belong to the $\spadesuit$-class. The proof for the other classes of function may be done straightforwardly in the same way.  
\begin{lem}\label{lem:pique}
Let $\cI$ be a familly of subsets of $\{1\/,\ldots\/,d\}$ and  $f_j^I$, $j\in \{1\/,\ldots\/, d\}$, $I\in \cI$,  be functions  either all increasing or all decreasing and they are either all positive or all negative then sums of product functions 
$$f(\bx) =\displaystyle \sum_{I\in \cI} \prod_{j\in I} f_j^I(x_j)$$
belong to the $\spadesuit$-class.
\end{lem}
\begin{proof}
 Let $A=\displaystyle\prod_{k=1}^d[a_k\/, b_k]$ and $j\in \{1\/,\ldots\/, d\}$ be fixed, for any $z\in [a_j\/,b_j]$,
 \begin{eqnarray*}
  \EE(f(z\/, \bX^{-j}) \ind_{\{\bX^ {-j} \in A^{-j}\}}) &= &\sum_{I\in \cI\/, \ j\in I} f_j^I(z) \EE\left(\prod_{k\in I\setminus \{j\}} f_k^I(X_k) \ind_{\{X_k \in [a_k\/,b_k]\}}\right) \\
  &&+ \sum_{I\in \cI\/, \ j\not\in I}\ \EE\left(\prod_{k\in I} f_k^I(X_k) \ind_{\{X_k \in [a_k\/,b_k]\}}\right) \/.
 \end{eqnarray*}
 So that $z\mapsto \EE(f(z\/, \bX^{-j}) \ind_{\{\bX^ {-j}\in A^{-j}\}})$ is constant on $[a_j\/,b_j]$ if and only if for any $z\/, z' \in [a_j\/,b_j]$,
 $$ \sum_{I\in \cI\/, \ j\in I} (f_j^I(z)-f_j^I(z')) \EE\left(\prod_{k\in I\setminus \{j\}} f_k^I(X_k) \ind_{\{X_k \in [a_k\/,b_k]\}}\right) =0 \/,$$
since the functions $f_k^I$ are   either all increasing or all decreasing and they are either all positive or all negative, this may happen  for all $z\/, z' \in [a_j\/, b_j]$, only if  $f_j^I$ is contant on $[a_j\/, b_j]$. 
\end{proof}

}

\end{appendix}

\end{document}